\documentclass[12pt]{amsart}

\usepackage{amssymb, amsmath, amscd,graphicx,subfigure}
\usepackage{amssymb,latexsym}

\usepackage{graphicx}

\theoremstyle{plain}

\newtheorem{lem}{Lemma}[section]

\newtheorem{prop}[lem]{Proposition}

\theoremstyle{definition}
\newtheorem{defi}[lem]{Definition}
\newtheorem{remn}[lem]{Remark}

\theoremstyle{plain}
\newtheorem{thm}{Theorem}
\newtheorem{cor}[thm]{Corollary}

\renewcommand{\int}{\operatorname{int}}

\newcommand{\cW}{\mathcal{W}}

\begin{document}

\title{Whitney towers and gropes in 4--manifolds}

\author{Rob Schneiderman}
\address{Courant Institute of Mathematical Sciences, New York University,
251 Mercer Street, New York NY 10012-1185 USA}
\email{schneiderman@courant.nyu.edu}
\thanks{The author is an NSF VIGRE postdoctoral fellow at the Courant
Institute of Mathematical Sciences.}






\begin{abstract}
Many open problems and important theorems in low-dimensional
topology have been formulated as statements about certain
2--complexes called {\em gropes}. This paper describes a precise
correspondence between embedded gropes in 4--manifolds and the
failure of the Whitney move in terms of iterated `towers' of
Whitney disks. The `flexibility' of these {\em Whitney towers} is
used to demonstrate some geometric consequences for knot and link
concordance connected to $n$-solvability, $k$-cobordism and grope
concordance. The key observation is that the essential structure
of gropes and Whitney towers can be described by embedded
unitrivalent trees which can be controlled during surgeries and
Whitney moves. It is shown that a Whitney move in a Whitney tower
induces an IHX (Jacobi) relation on the embedded trees.

\end{abstract}

\maketitle

\section{Introduction}
Many open problems and important theorems in low-dimensional
topology, including the classification theory of topological
4--dimensional manifolds and the study of knots and links in
3--manifolds, have been formulated in terms of statements about
smooth maps of certain 2--complexes called {\em gropes}. As
suggested by its name, a grope is built inductively from layers of
surfaces which ``reach into'' a 4--manifold in an attempt to
approximate an embedded 2--disk (see \cite{T1,T2}). The existence of
generically embedded 2--disks in dimensions greater than 4 allows
for use of the Whitney move, a vital part of the surgery programs
which yield classification theorems for higher dimensional
manifolds and knotted submanifolds. The failure in general of the
Whitney move in dimension 4 is a defining characteristic of
low-dimensional topology and this paper describes a precise
correspondence between gropes and an approximation to the Whitney
move via certain ``towers'' of iterated Whitney disks.

\begin{figure}[ht!]
        \centerline{\includegraphics[scale=.5]{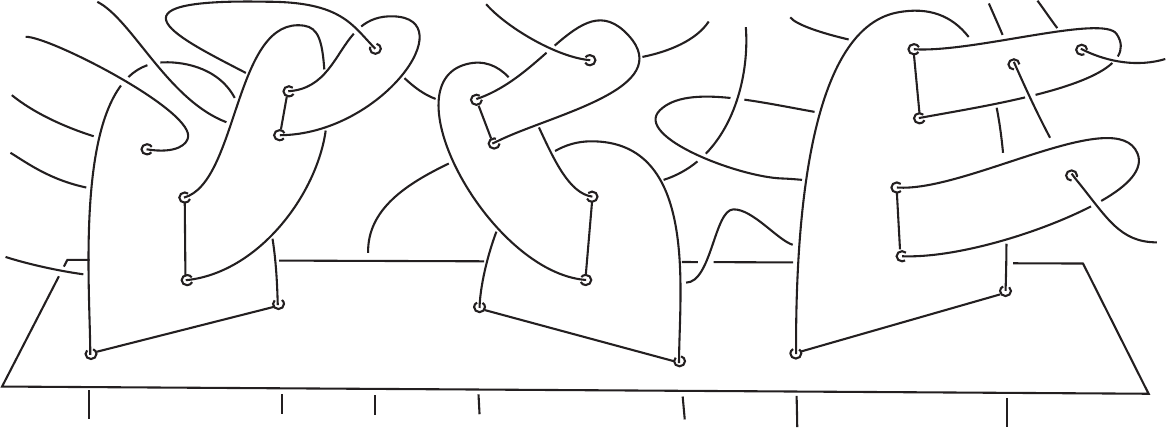}}
        \caption{Part of a Whitney tower.}
        \label{W-tower5-fig}

\end{figure}

A Whitney move eliminates a pair of singularities between immersed surfaces in a 4--manifold but
will also create new intersections if the guiding Whitney disk contains singularities in its
interior. A {\em Whitney tower} is constructed by pairing such singularities in the interiors of
Whitney disks with ``higher order'' Whitney disks in an attempt to find embedded Whitney disks via
``higher order'' Whitney moves (Figure~\ref{W-tower5-fig} and Section~\ref{w-tower-sec}).

The {\em order} of a Whitney tower (Definition~\ref{w-tower-defi}) is determined by the number of
layers of Whitney disks added to the immersed surfaces and the {\em class} of a grope
(Figure~\ref{class4gropeA-fig} and Definition~\ref{grope-defi}) measures its complexity in terms of
the number of its layers or {\em stages} of surfaces.

In the case of gropes, all singularities are usually contained in {\em caps} which are 2--disks
that are mapped in after all stages of embedded surfaces have been attached. The caps are attached
along essential curves called {\em tips} (details in Section~\ref{grope-sec}).

The correspondence between class and order is described by the following basic version of our main
result:

\begin{thm}\label{class-order-thm}
For any collection of embedded closed curves $\gamma_i$ in the boundary of a 4--manifold $X$, the
following are equivalent:
\begin{enumerate}

\item The $\gamma_i$ bound disjoint properly embedded class $n$
gropes $g_i$ with null-homotopic tips in $X$.

\item The $\gamma_i$ bound properly immersed 2--disks $D_i$ admitting an order $(n-1)$
Whitney tower $\cW$ in $X$.

\end{enumerate}

\end{thm}
The proof of Theorem~\ref{class-order-thm} is given by a local
construction meaning that, given caps for the $g_i$, the $D_i$ and
$\cW$ are constructed in a neighborhood of the capped gropes
$g^c_i$ and, given the $D_i$ contained in $\cW$, the gropes $g_i$
(and caps) are constructed in a neighborhood of $\cW$. (This
construction also applies to {\em closed} (boundaryless) gropes
and Whitney towers.)
\subsection*{Unitrivalent trees}
In fact, a much more precise and more general correspondence between gropes and Whitney towers can
be described in terms of unitrivalent trees (\ref{grope-trees}, \ref{int-point-trees}) which
capture the essential structure of both gropes and Whitney towers. This correspondence is given in
Theorem~\ref{tower-to-grope-thm} and Theorem~\ref{grope-to-tower-thm}, which are stated and proved
in Section~\ref{sec:thms} (and directly imply Theorem~\ref{class-order-thm}). The key element used
in the proofs of these theorems and their corollaries is that the associated trees, embedded as
subsets of the gropes and Whitney towers, can be preserved and controlled during modifications such
as ambient surgery or Whitney moves.

\begin{remn}\label{rem:tau}
Besides providing an interpretation of gropes in terms of Whitney
disks, the notion of Whitney towers comes equipped with an
obstruction theory that generalizes the usual intersection theory
for surfaces in 4--manifolds. The associated invariants $\tau_n$
take values in graded groups of trees which generalize the
well-known target groups of finite type concordance invariants of
links. These higher order intersection trees are enhancements of
the trees discussed in this paper and are used in \cite{ST3} to
give a geometric description of Milnor's $\mu$--invariants and the
(reduced) Kontsevich integral for links, and are also used in
\cite{ST4} to obtain results on separating homotopy classes of
surfaces in 4--manifolds. See also the more recent papers \cite{CST2,CST5,WTCCL}.
\end{remn}

\subsection*{Symmetric gropes and $\mathbf{n}$--solvability}
It should be noted that the Whitney towers described here differ from the towers that appear in the
theory of Casson handles \cite{F,Q} in that the latter always involve killing the
so-called ``accessory circles'' that run through the singularities. The gropes that arise in the
context of Casson handles are certain capped {\em symmetric} gropes whose singularities are
restricted to intersections among the caps \cite{FQ}.

The notion of Whitney tower described here {\em is} related to
general symmetric gropes (without restrictions on caps) and the
recently discovered filtration of the classical knot concordance
group in terms of {\em n-solvability} by Cochran, Orr and Teichner
\cite{COT}: A Whitney tower of {\em height} $n$ (resp. $n.5$) as
defined in \cite{COT} is a Whitney tower of order $2^n-2$ (resp.
$2^n+2^{(n-1)}-2$) with certain additional restrictions on the
type of allowable intersections (see
Section~\ref{sec:height-cor}). This definition of height is in
rough correspondence with the usual definition of height for a
symmetric (uncapped) grope and in \cite{COT} it is shown that a
knot in the 3--sphere $S^3=\partial B^4$ is $n$-solvable if it
bounds an embedded symmetric grope of height $n+2$ or a Whitney
tower of height $n+2$ in $B^4$. It is not known if these geometric
conditions are equivalent to each other (or to being $n$-solvable)
but we do have:
\begin{cor}\label{cor:height}
If a knot in $S^3$ bounds a properly embedded grope of height $n$ (resp. $n.5$) in $B^4$, then the
knot bounds a Whitney tower of height $n$ (resp. $n.5$) in $B^4$.
\end{cor}
The proof of Corollary~\ref{cor:height} suggests that bounding a grope of height $n$ is probably a
stronger condition than bounding a Whitney tower of height $n$ (see Section~\ref{sec:height-cor}).


\subsection*{Half-gropes, $\mathbf{k}$--cobordism and geometric IHX}
In the general study of knot and link concordance (cobordism) the translation between gropes and
Whitney towers can be helpful in understanding connections between algebraically and geometrically
defined filtrations. This is the case in the next two corollaries which are proved by applying a
geometric ``IHX'' Jacobi relation to manipulate the trees associated to a Whitney tower
(Lemma~\ref{IHX-lemma}).

In analogy with the fact that the $n$th term of the lower central series of a group is generated by
{\em simple} (right- or left-normed) commutators \cite{MKS} we have the following geometric
result which implies in particular that class $n$ {\em grope concordance} is generated by class $n$
(annulus-like) {\em half-gropes} (details in Section~\ref{sec:half-grope-cor}):
\begin{cor}\label{cor:half-grope}
If $\gamma_i$ bound disjoint properly embedded $A_i$--like class $n$ gropes $g_i$ with
null-homotopic tips in a 4--manifold $X$, then the $\gamma_i$ bound disjoint properly embedded
$A_i$--like class $n$ {\em half-gropes} $h_i$ with null-homotopic tips in $X$, with the $h_i^c$
contained in a neighborhood of $g_i^c$ for any choice of caps on the $g_i$.
\end{cor}
Since a symmetric grope of height $n$ has class $2^n$, Corollary~\ref{cor:half-grope} also gives a
geometric (embedded) analogue of the fact that the $n$th derived subgroup is contained in the
$2^n$th lower central subgroup of a group.

In the case of knots in $S^3$, it was shown constructively in \cite{S2} (see also \cite{CT2}) that
the Arf invariant is the only obstruction to the existence of a (half)-grope concordance of
arbitrarily high class. Thus, the Von Neumann signatures of \cite{COT} which obstruct
$n$-solvability are obstructions to ``inverting'' the geometric manipulations of gropes and Whitney
towers via geometric IHX constructions (as used, for instance, in the proof of
Corollary~\ref{cor:half-grope} above) to convert class into height.

The filtration of link concordance classes by the notion of $k$-{\em cobordism} was introduced and
studied in \cite{C1}, \cite{C2}, \cite{C3} and \cite{O1}. In particular, a link $L$ in $S^3$ is
$k$-{\em slice} ($k$-null-cobordant) if the link components bound disjoint surfaces in $B^4$ which
``look (algebraically) like'' slice disks modulo the $k$th term of the lower central series of the
link group $\pi_1(S^3- L)$. This was shown in \cite{IO} to be equivalent to $L$ having vanishing
Milnor $\overline{\mu}$-invariants up through length $2k$. The precise relation between grope
concordance and Milnor's invariants is not known, however class $n$ grope concordance implies
length $n$ Milnor--equivalence \cite{KrT}. Evidence that class $n$ grope concordance is {\em
stronger} than length $n$ Milnor--equivalence is provided by the easy proof of the following
corollary, in contrast to the difficulties encountered in the just mentioned result in \cite{IO}.

\begin{cor}\label{cor:k-slice}
If the components of a link $L$ in $S^3$ bound disjoint properly embedded class $2k$ gropes in the
$B^4$, then $L$ is $k$-slice.
\end{cor}
The proof of Corollary~\ref{cor:k-slice} uses the flexibility of Whitney towers to ``evenly
distribute'' the gropes' higher surface stages over symplectic sets of circles on the bottom
surfaces. In fact, the proof yields a stronger conclusion as described in
Section~\ref{sec:k-slice}. (The same proof also shows more generally that class $2k$ grope
concordance implies $k$-cobordism.)

\subsection*{Further applications}
Although closely related, Whitney towers are slightly more general
objects than gropes as is suggested by the association of {\em
rooted} trees to gropes and {\em un}rooted trees to Whitney
towers. This generality is conducive to defining geometric
invariants associated to an obstruction theory (as mentioned in
Remark~\ref{rem:tau}) as well as manipulating the shape of gropes
(as in Corollaries~\ref{cor:half-grope} and \ref{cor:k-slice}
above). However, in certain constructions it is useful to convert
Whitney towers into gropes, for instance to take advantage of nice
properties of embedded grope complements. This interplay is
exploited in \cite{CST2,WTCCL} which describes the geometry of Milnor's
$\mu$--invariants in terms of Whitney towers and grope
concordance.
 Recent work of Conant and Teichner \cite{CT1,CT2} suggests that the study of both
3-- and 4--dimensional grope cobordism of restricted graph type can play an important role in the
general theory of knots and links. The relevant graphs in 4--dimensions are the unitrivalent trees
which occur here naturally in the context of Whitney towers. Perhaps an analogous notion of Whitney
tower projected into 3-dimensions with associated unitrivalent graphs could be useful in
understanding 3--dimensional grope cobordism.

\subsection*{Outline}
Gropes are discussed in Section~\ref{grope-sec}. Section~\ref{w-tower-sec} covers Whitney towers,
including two lemmas (\ref{split-tower-lem} and \ref{subtower-lemma}) which are used in subsequent
constructions. Section~\ref{sec:grope-subtowers} introduces the hybrid ({\em split}) {\em grope
subtowers} which are used to interpolate between gropes and Whitney towers.
Theorem~\ref{class-order-thm} is then proved via the more detailed Theorem~\ref{tower-to-grope-thm}
and Theorem~\ref{grope-to-tower-thm} which are stated and proved in Section~\ref{sec:thms}. A proof
of Corollary~\ref{cor:height} is given in Section~\ref{sec:height-cor}. A proof of
Corollary~\ref{cor:half-grope} is given in Section~\ref{sec:half-grope-cor}, which also contains a
geometric (IHX) construction (Lemma~\ref{IHX-lemma}) which is used in \cite{CST2,ST3}. A proof of
Corollary~\ref{cor:k-slice} is given in Section~\ref{sec:k-slice}.

\subsection*{Conventions}\label{conventions}
All maps and manifolds are assumed smooth and oriented. Surfaces in 4--manifolds are illustrated in
figures showing a 3--dimensional slice of local 4--dimensional coordinates with the understanding
that sheets of surfaces that appear as 1-dimensional arcs in the ``present'' 3--dimensional slice
extend as a product into the ``past and future'' coordinate. Sheets of surfaces which are contained
in the 3--dimensional slice may appear either translucent or opaque.
\begin{figure}[ht!]
        \centerline{\includegraphics[scale=.5]{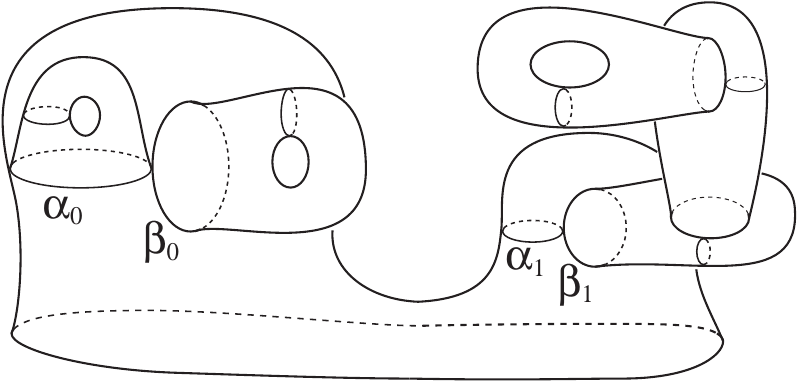}}
        \caption{A grope of class 4 (See Definition~\ref{grope-defi}).}
        \label{class4gropeA-fig}

\end{figure}
\section{Gropes}\label{grope-sec}
This section contains mostly standard grope terminology for the purpose of fixing notation since
there are subtle differences in definitions and approaches throughout the literature. See any of
\cite{FQ}, \cite{FT} or \cite{KrT} for detailed discussions of gropes in 4--manifolds.

There will be however be a few mildly non-standard wrinkles which are worth pointing out: It will
be convenient to consider the punctured {\em starting surface} $A^0$ of an $A$--like grope $g$ as a
{\em 0th stage} of $g$,  and it will be helpful think of $g$ as formed from $A^0$ by attaching many
genus one gropes rather than a single higher genus grope. In light of these notational conventions,
it turns out to be convenient to associate to each grope a disjoint union of unitrivalent trees (as
in \cite{CST}, \cite{CT1} and \cite{CT2}) rather than the customary (single) multivalent tree (e.g.
in \cite{KrT}). This point of view is in line with Krushkal's grope-splitting technique \cite{Kr}
which simplifies the combinatorics of gropes. Note also that what we will refer to as a {\em
dyadic} grope is called a ``grope with dyadic branches'' in \cite{KQ}.
\begin{figure}[ht!]
        \centerline{\includegraphics[scale=.5]{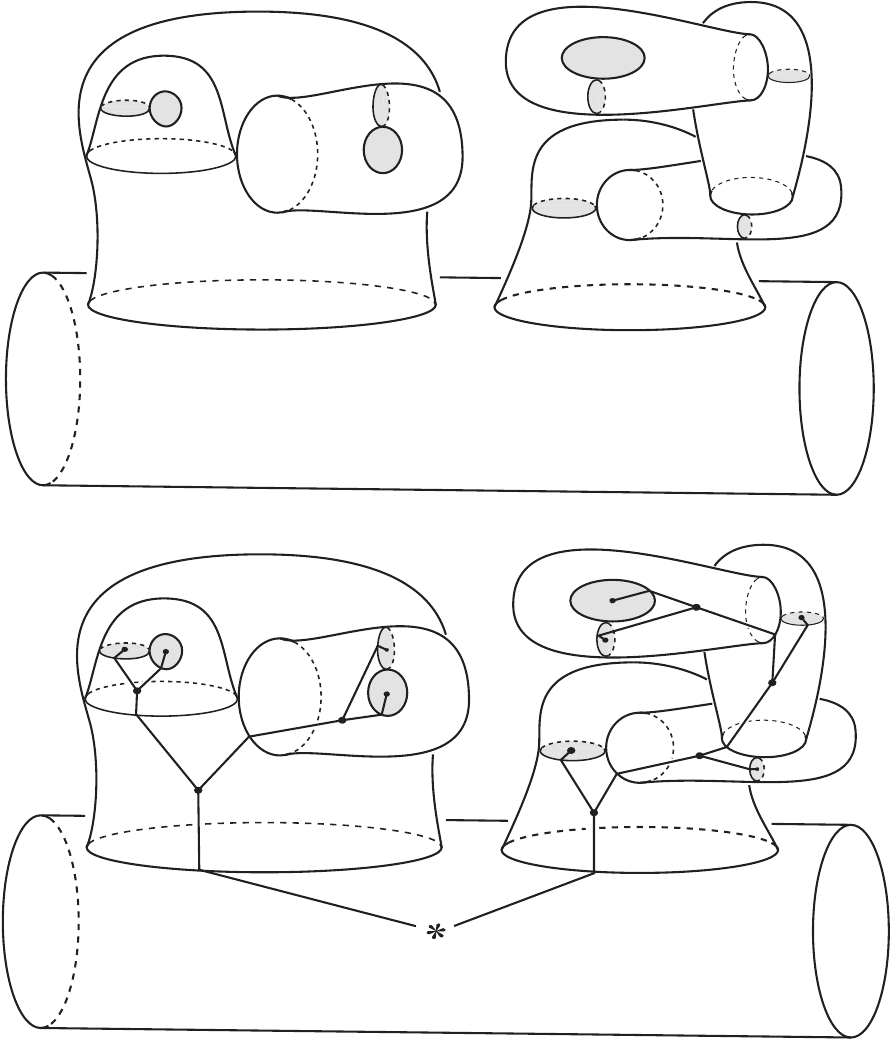}}
        \caption{A dyadic annulus-like capped grope $g^c$ of class 4, also shown
        with the embedding of $t(g^c)$ in $g^c$.
        The 0th stage of $g^c$ is the twice-punctured annulus.}
        \label{capped-A-like-grope1-and-with-trees-fig}

\end{figure}
\subsection{Grope terminology}\label{grope-terminology}

\begin{defi}[\cite{FT}]\label{grope-defi}
A {\em grope} is a special pair (2--complex, circle). A grope has a {\em class}
$n\in\{1,2,\ldots,\infty\}$. A class 1 grope is defined to be the pair (circle, circle). A class 2
grope $(S,\partial S)$ is a compact oriented connected surface $S$ with a single boundary
component. For $n>2$, a grope of class $n$ is defined inductively as follows: Let
$\{\alpha_j,\beta_j, j=1,\ldots,\mbox{genus}\}$ be a chosen standard symplectic basis of circles
for a class~2 grope $S$. For any positive integers $a_j,b_j$, with $a_i+b_j\geq n$ and
$a_{j_0}+b_{j_0}=n$ for at least one index $j_0$, a grope of class $n$ is formed by attaching a
class $a_j$ grope to each $\alpha_j$ and a class $b_j$ grope to each $\beta_j$ (See
Figure~\ref{class4gropeA-fig}).
\end{defi}
Here ``attaching a class 1 grope'' is understood to mean ``not attaching a grope at all''. The
surfaces in a grope $g$ are called {\em stages} and the basis circles in (all stages of) $g$ which
do not have a surface stage attached to them are the {\em tips} of $g$. Attaching 2--disks, called
{\em caps}, to all the tips of $g$ yields a {\em capped grope} $g^c$ and the underlying uncapped
grope $g$ is the {\em body} of $g^c$.

It is customary to omit the boundary of a grope from notation. We adopt the convention that the
{\em tip} of a class 1 grope is the grope (a circle) itself, so a class one capped grope is just a
circle bounding a disk.

Note that a class $n$ grope is {\em not} really also a grope of class $m$ for $m<n$, but can always
be made into one by deleting some surface stages (in order to satisfy $a_{j_0}+b_{j_0}=m$ for at
least one index $j_0$).

\subsection{Dyadic {\boldmath $A$}--like gropes with 0th
stages.}\label{dyadic-gropes-defi} The gropes of Definition~\ref{grope-defi} are special cases of
the more general {\em $A$--like (capped) grope} of class $n$ which is gotten from a starting
surface $A$ by replacing disks in $A$ with (capped) gropes of class $n$ as defined above. All the
notions such as stages, tips and caps apply to $A$--like gropes without change with the following
addition: If $g$ is an $A$--like grope then the {\em 0th stage} of $g$ is defined to be the
punctured surface $A^0$ which is $A$ minus the disks that were replaced by gropes to form $g$. The
0th stage is included in the body of an $A$--like capped grope. In this language, the gropes of
Definition~\ref{grope-defi} can be thought of as disk-like gropes.

An $A$--like (capped) grope $g$ is {\em dyadic} if all surface stages of all gropes attached to the
0th stage $A^0$ are genus one surfaces (See Figure~\ref{capped-A-like-grope1-and-with-trees-fig}).
Since we are allowing $A^0$ to have many punctures, any (capped) $A$--like grope can be converted
to a (capped) dyadic $A$--like grope by Krushkal's grope-splitting technique \cite{Kr}.

The notion of dyadic gropes was introduced in \cite{KQ}; our terminology is slightly different in
that instead of allowing higher genus in the first stage we attach many genus 1 gropes to the 0th
stage. Also, for brevity we are using the term ``dyadic grope'' here instead of (the more precise)
``grope with dyadic branches'' as in \cite{KQ}.

\subsection{Rooted trees for gropes}\label{grope-trees}
A {\em tree} is a connected graph without 1--cycles. A {\em rooted} tree has a single preferred
univalent vertex called the {\em root}. To each $A$--like dyadic capped grope $g^c$, we associate a
disjoint union $t(g^c)$ of rooted unitrivalent trees which is essentially the dual one-complex: The
vertices of $t(g^c)$ are (dual to) the stages and caps of $g^c$, with two vertices joined by an
edge if the corresponding stages/caps meet in a circle. The univalent vertices corresponding to the
0th stages are the roots.

It will be helpful to think of $t(g^c)$ as being embedded in $g^c$ in the following way. Choose a
basepoint in each stage (including $A^0$) and a basepoint in each cap of $g^c$. If $g^c$ was formed
by removing just one disk from $A$, then connecting basepoints in adjacent stages and caps by
sheet-changing paths yields an embedded connected unitrivalent tree; here ``adjacent'' means
``intersecting along a circle'' so that dual stages/caps (having a single intersection point in
their boundaries) are not considered to be adjacent. If $g^c$ was formed by replacing $m$ disks of
$A$, then we get $m$ unitrivalent trees (each ``sprouting from'' the basepoint in the $0$th stage
$A^0$) and $t(g^c)$ is the disjoint union of these $m$ trees
(Figure~\ref{capped-A-like-grope1-and-with-trees-fig}). Each trivalent vertex of $t(g^c)$
corresponds to a genus one stage in $g^c$ and each univalent vertex of $t(g^c)$ corresponds to a
cap  of $g^c$, except for one univalent vertex on each connected component of $t(g^c)$ which
corresponds to $A^0$. These univalent vertices on the $A^0$ are the root vertices of $t(g^c)$. To
the underlying uncapped grope $g$ is associated the same disjoint union of rooted trees denoted
$t(g)$.

Following the terminology of \cite{ST3}, define the {\em order} of
a unitrivalent tree to be the number of trivalent vertices. Note
that the class of $g^c$ is equal to one more than the minimum of
the orders of the trees in $t(g^c)$.

The reader familiar with previous associations of multi-valent trees to arbitrary $A$--like gropes
(e.g. \cite{KrT}) can check that our definition of $t(g^c)$ is essentially what you would get after
applying Krushkal's grope-splitting procedure. In fact, a formal sum of {\em vertex-oriented} trees
can be associated to an {\em oriented} grope \cite{CST,CT1,CT2}; although we will
not work with orientations in this paper, the notation here has been chosen to be compatible with
the just cited works.

\subsection{Proper immersions of gropes}\label{immersions}
A surface is {\em properly} immersed in a 4--manifold if boundary is embedded in boundary and
interior is immersed in interior. Immersions of (capped) gropes into a 4--manifold are required to
factor through an embedding in 3--space followed by standard product thickenings and plumbings, so
that a regular neighborhood of the immersion contains disjoint parallel copies of any embedded
subsets of the (capped) grope. An immersion of (capped) gropes is {\em proper} if the bottom stage
surfaces are properly immersed and all other stages (and caps) are immersed in the interior.

In much of the literature, a proper immersion of a capped grope in a 4--manifold will have
restrictions on allowable grope/cap singularities. We will not make such restrictions, however it
will be convenient to arrange for all intersections to occur among $0$th stages and between caps
and $0$th stages .

\section{Whitney towers}\label{w-tower-sec}
After defining {\em Whitney towers}, the main goal of this section
is to show how the essential geometric structure of a Whitney
tower $\cW$ is captured by a disjoint union $t(\cW)$ of
unitrivalent (labelled) trees. The surfaces of $\cW$ are indexed
by brackets which correspond to rooted trees, and each
intersection point $p$ of $\cW$ is then assigned an unrooted tree
$t(p)$ which is a pairing of the rooted trees that correspond to
the intersecting surfaces. After introducing the notion of ({\em
split}) {\em subtowers}, Lemma~\ref{split-tower-lem} then shows
how $\cW$ can be decomposed into essential parts which are
described by the disjoint union $t(\cW)$ of all the $t(p)$.
Finally, Lemma~\ref{subtower-lemma}, which will play a key role in
later proofs, gives further evidence of the essential nature of
$t(\cW)$ by showing how the trees $t(p)$ are ``stable'' under
certain Whitney moves.

All of the material in this section can be enhanced to take into account orientations as well as
the fundamental group of the ambient manifold by adding vertex-orientations and edge-decorations to
the trees. Such enhancements are used in the obstruction theory of \cite{ST1}, \cite{ST3},
\cite{ST3}, and \cite{WTCCL}. Notation in this paper has been chosen to be consistent with these
papers where they overlap.
\begin{figure}[ht!]
        \centerline{\includegraphics[scale=.5]{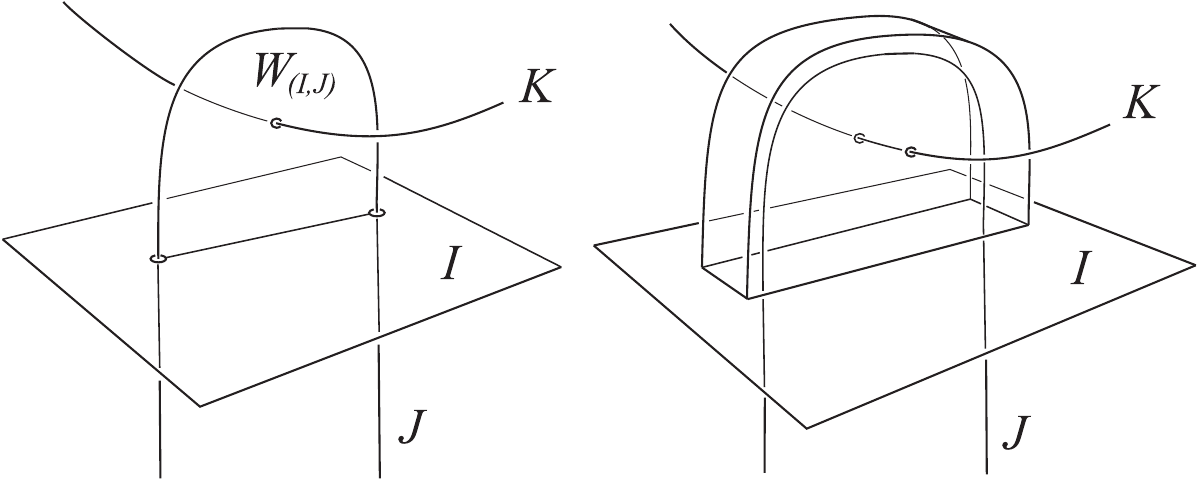}}
        \caption{The embedded Whitney disk $W_{(I,J)}$ guides a Whitney move that
        eliminates a cancelling
        pair of intersections between surface sheets $I$ and $J$ via an isotopy of $I$.
        Since $W_{(I,J)}$ had an interior intersection with a surface sheet $K$,
        this Whitney move creates a new cancelling
        pair of intersections between $I$ and $K$.}
        \label{W-disk-move2-fig}

\end{figure}
\subsection{Whitney disks}\label{W-disk-subsec}
We refer the reader to \cite{FQ} for a detailed description of the {\em Whitney move} in dimension
4. For our purposes it is enough to understand the model Whitney move on an embedded Whitney disk
in 4--space as illustrated in Figure~\ref{W-disk-move2-fig}, which shows the effect of a Whitney
move on three local sheets of surfaces: The Whitney move is guided by an embedded {\em Whitney
disk} labelled $W_{(I,J)}$ which {\em pairs} a (geometrically) {\em cancelling pair} of
intersection points between sheets of surfaces labelled $I$ and $J$. In
Figure~\ref{W-disk-move2-fig}, the $I$-sheet is locally contained in the 3--dimensional ``present''
and the $J$- and $K$-sheets extend into ``past and future''. The Whitney move eliminates the
cancelling pair of intersections between $I$ and $J$ by an isotopy of $I$ across $W_{(I,J)}$ at the
cost of introducing a new cancelling pair of intersection points between $I$ and the sheet labelled
$K$ which intersected the interior of the guiding Whitney disk $W_{(I,J)}$.

\subsection{Whitney towers}\label{W-tower-sub-sec}
In general, a Whitney disk may have multiple interior self-intersections and intersections with
other surface sheets; however we will require that arbitrary Whitney disks resemble the model near
their boundaries. By pairing up interior intersections with higher order Whitney disks we are led
to the notion of a Whitney tower:
\begin{defi}\label{w-tower-defi}\mbox{}
\begin{itemize}
\item A {\em surface of order 0} in a 4--manifold $X$
is a properly immersed surface (boundary embedded in the boundary
of $X$ and interior immersed in the interior of $X$). A {\em
Whitney tower of order 0} in $X$ is a collection of order 0
surfaces.
\item The {\em order of a (transverse) intersection point} between a surface of order $n$ and a
surface of order $m$ is $n+m,$.
\item The {\em order of a Whitney disk} is $(n+1)$ if it pairs intersection points of order $n$.
\item For $n\geq 0$,
a {\em Whitney tower of order $(n+1)$} is a Whitney tower $\cW$ of
order $n$ together with Whitney disks pairing all order $n$
intersection points of $\cW$. The interiors of these top order
disks are allowed to intersect each other as well as lower order
surfaces.
\end{itemize}
The Whitney disks in a Whitney tower are required to be {\em
framed} (see \cite{FQ}) and have disjointly embedded boundaries.
It will also be assumed that the order 0 surfaces are {\em
0-framed} (see 1.2 of \cite{FQ}).
\end{defi}
Thus, in an order $n$ Whitney tower all intersection points of
order less than $n$ occur in cancelling pairs with respect to
(arbitrarily) chosen orientations of all Whitney disks (see
Figure~\ref{W-tower5-fig}). Note that the boundary of any Whitney
disk is not allowed to change sheets except at the intersection
points paired by the Whitney disk.

Some further terminology: If $\cW_0$ is an order 0 Whitney tower
and there exists an order $n$ Whitney tower $\cW_n$ containing
$\cW_0$ as its order 0 surfaces, then $\cW_0$ is said to {\em
admit} an order $n$ Whitney tower and any one of the order 0
surfaces in $\cW_n$ is said to {\em support} the $n$th order
Whitney tower $\cW_n$.
\begin{figure}[ht!]
        \centerline{\includegraphics[scale=.5]{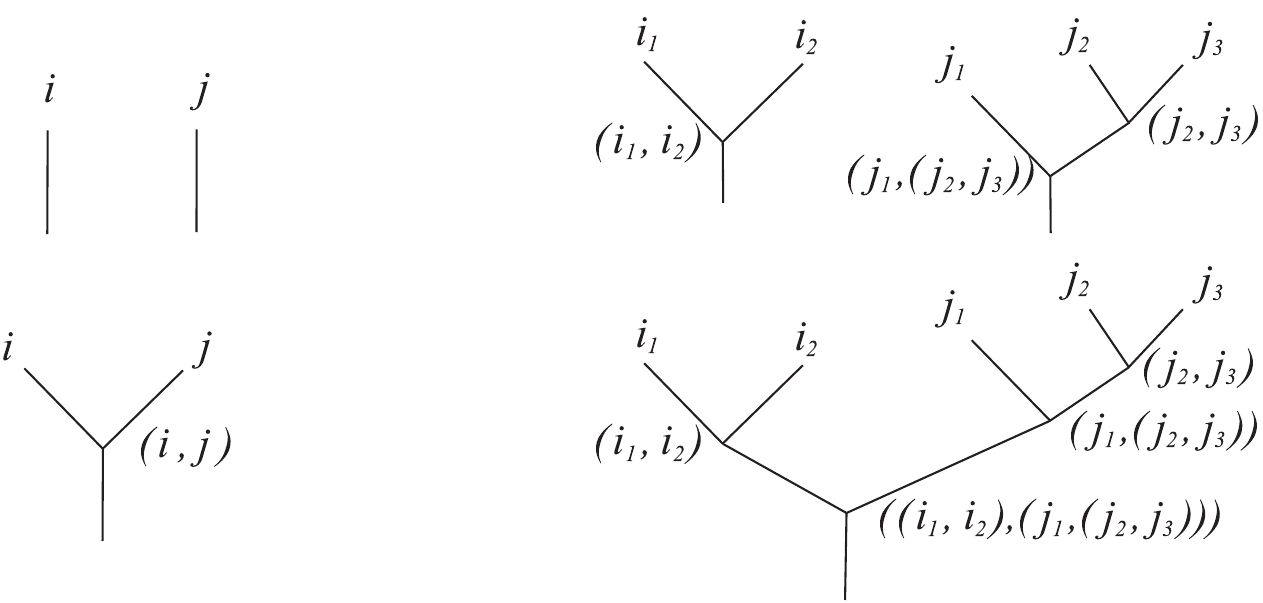}}
        \caption{Rooted trees and brackets: $t(i)$ and $t(j)$ (upper left) and their rooted product
        $t(i,j)=t(i)\ast t(j)$ (lower left);
        $t(i_1,i_2)$ and $t(j_1,(j_2,j_3))$ (upper right) and their rooted product
        $t((i_1,i_2),(j_1,(j_2,j_3)))=t(i_1,i_2)\ast t(j_1,(j_2,j_3))$
        (lower right).}
        \label{bracket-treesA-fig}

\end{figure}
\subsection{Rooted trees and brackets}\label{trees-brackets}
Non-associative but commutative (unordered) {\em bracketings} of
elements from some index set correspond to rooted labelled
(unoriented) unitrivalent trees as follows. A bracketing $(i)$ of
a singleton element $i$ from the index set corresponds to the
rooted {\em chord} $t(i)$ having a single edge with one vertex
labelled by $i$ and the other vertex designated as the root. A
bracketing $(I,J)$ of brackets $I$ and $J$ corresponds to the {\em
rooted product} $t(I) \ast t(J)$ of the trees $t(I)$ and $t(J)$
which identifies together the roots of $t(I)$ and $t(J)$ to a
single vertex and ``sprouts'' a new rooted edge at this vertex
(Figure~\ref{bracket-treesA-fig}). Thus, the (non-root) univalent
vertices of the tree $t(I)$ associated to a bracket $I$ are
labelled by elements from the index set and the trivalent vertices
correspond to sub-bracketings of $I$, with the trivalent vertex
adjacent to the root corresponding to $I$.
\begin{figure}[ht!]
        \centerline{\includegraphics[scale=.50]{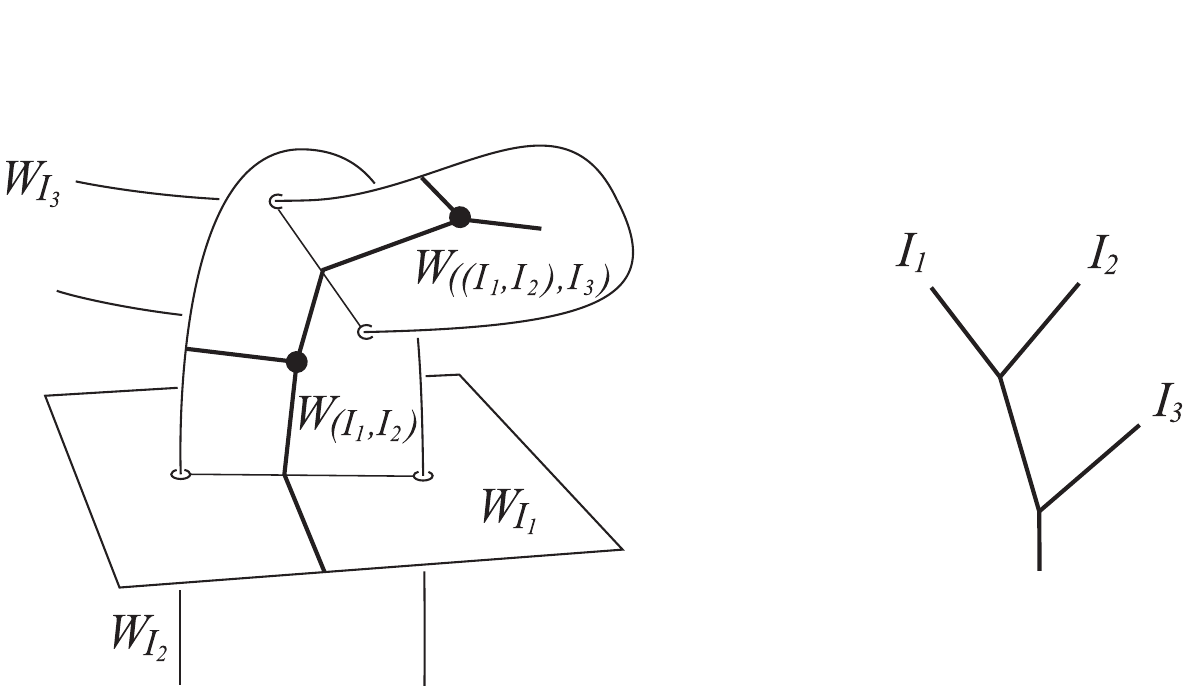}}
        \caption{A Whitney disk $W_{((I_1,I_2),I_3)}$ and (part of) its
        associated tree $t(W_{((I_1,I_2),I_3)})$ shown both as a subset of
        the Whitney tower and as an abstract rooted tree. (If any $I_i$ is not
        a singleton, then that label $I_i$ represents a sub-tree.)}
        \label{W-disk-treeC-fig}

\end{figure}
\subsection{Rooted trees for Whitney disks}\label{w-disk-tree-subsec}
The rooted trees and bracketings of the previous subsection are
associated to the surfaces in a Whitney tower $\cW$ both formally
and geometrically as follows: A bracketing $(i)$ of a singleton
element $i$ from the index set is associated to each order zero
surface. The bracket $(I,J)$ is associated to a Whitney disk
pairing intersections between surfaces with associated brackets
$I$ and $J$. Using brackets as subscripts, we write $A_i$ for an
order zero surface (dropping the brackets around the singleton
$i$) and $W_{(i,j)}$ for a first order Whitney disk that pairs
intersections between $A_i$ and $A_j$. In general, we write
$W_{(I,J)}$ for a Whitney disk pairing intersections between $W_I$
and $W_J$, with the understanding that if a bracket $I$ is just a
singleton $(i)$, then the surface $W_I=W_{(i)}$ is just the order
zero surface $A_i$. The rooted labelled tree $t(W_I)$ associated
to $W_I$ is defined to be $t(I)$, the tree that corresponds to the
bracket $I$ as before (\ref{trees-brackets}). Note that the order
of $W_I$ is equal to the order of $t(W_I)$ (i.e. the number of
trivalent vertices).

In fact, the tree $t(W_I)$ can be mapped into $\cW$ in a similar
fashion to the case of gropes (\ref{grope-trees}): First fix a
basepoint in the interior of each surface (including the Whitney
disks) of $\cW$. (If an order zero surface is not connected put a
basepoint in each component.) Now map the vertices (other than the
root) of $t(W_I)$ to the basepoints of the surfaces corresponding
to the sub-brackets of $I$ and map the edges (other than the edge
adjacent to the root) of $t(W_I)$ to sheet-changing paths between
these basepoints, as illustrated in Figure~\ref{W-disk-treeC-fig}.
Then embed the root and its edge anywhere in the interior of
$W_I$.
\begin{figure}[ht!]
        \centerline{\includegraphics[scale=.45]{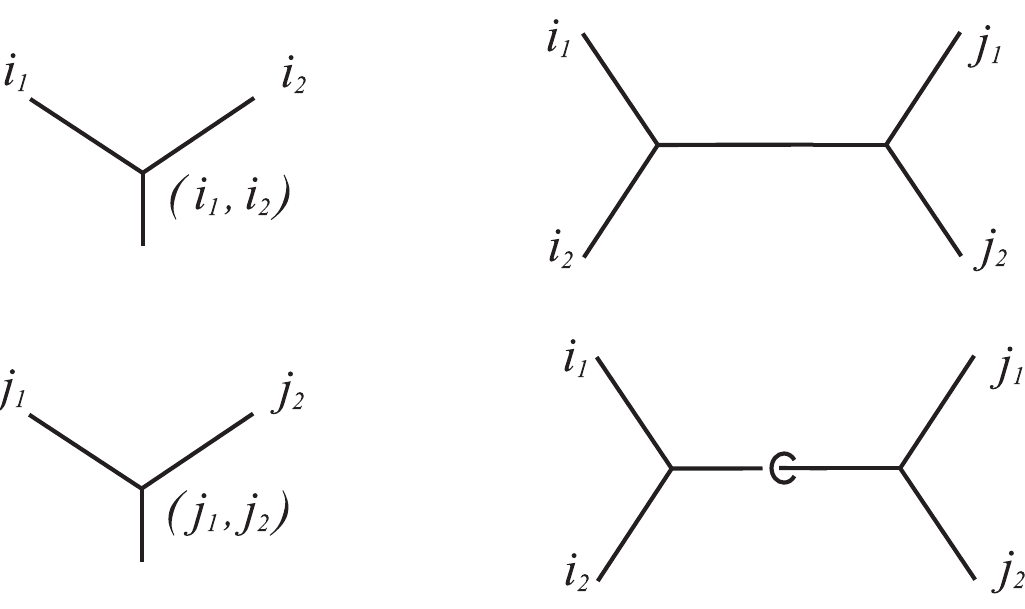}}
        \caption{On the left, a pair of rooted trees $t(I)$ and $t(J)$
        corresponding to first order Whitney disks $W_{I}$ and
        $W_{J}$ with $I=(i_1,i_2)$ and $J=(j_1,j_2)$. On the upper
        right, the inner product $t(p)=t(I)\cdot t(J)$ associated to
        a 2nd order intersection point $p\in W_{I}\cap W_{J}$ and on the
        lower right, the punctured tree $t^{\circ}(p)$ that also
        keeps track of $p$.}
        \label{inner-product-trees-fig}

\end{figure}
\subsection{Trees for intersection points}\label{int-point-trees}

Given a pair $t(I)$ and $t(J)$ of rooted trees, define the {\em
inner product} $t(I)\cdot t(J)$ to be the labelled {\em un}rooted
tree gotten by identifying together the root vertices of $t(I)$
and $t(J)$ to a single (non-vertex) point. The tree $t(p)$
associated to a (transverse) intersection point $p\in W_I\cap W_J$
between surfaces $W_I$ and $W_J$ in a Whitney tower $\cW$ is
defined to be the inner product $t(W_I)\cdot t(W_J)$ ($=t(I)\cdot
t(J)$) of the rooted trees corresponding to $W_I$ and $W_J$ as
illustrated in Figure~\ref{inner-product-trees-fig}. Note that the
order of $p$ is equal to the order of $t(p)$.

\begin{figure}[ht!]
        \centerline{\includegraphics[scale=.5]{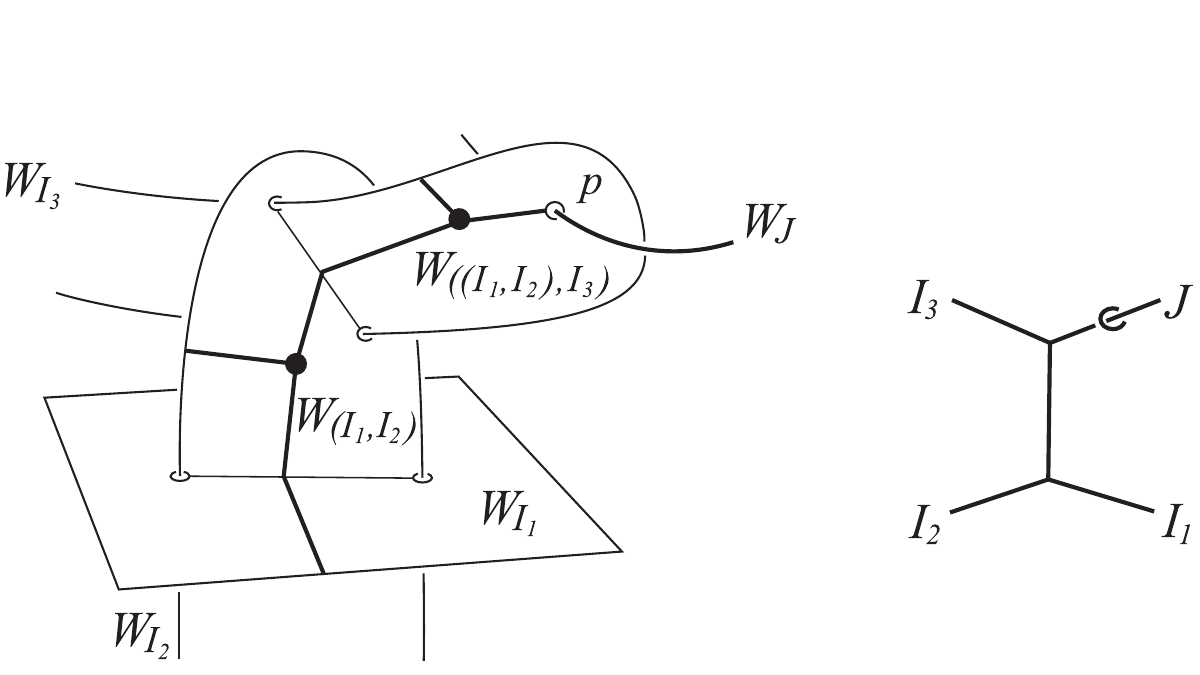}}
        \caption{A local picture of the punctured tree $t^{\circ}(p)$ associated to
        an intersection point $p\in W_I\cap W_J$, shown as a subset
        of the Whitney tower and as an abstract labelled (punctured) tree.}
        \label{W-disk-int-point-treeC-fig}

\end{figure}

The above mentioned mappings of $t(W_I)$ and $t(W_J)$ in $\cW$
give rise to a mapping of $t(p)$ into $\cW$: Just map the root
vertices of $W_I$ and $W_J$ to $p$ and the adjacent edges become a
path between the basepoints of $W_I$ and $W_J$ which changes
sheets at $p$ (Figure~\ref{W-disk-int-point-treeC-fig}). This
mapping can be taken to be an embedding of $t(p)$ into $\cW$ if
all the Whitney disks ``beneath'' $W_I$ and $W_J$ (corresponding
to sub-brackets of $I$ and $J$) are distinct.

It is sometimes convenient to keep track of the edge of $t(p)$
that corresponds to $p$ by marking that edge with a small linking
circle as in Figure~\ref{inner-product-trees-fig} and
Figure~\ref{W-disk-int-point-treeC-fig}; such a {\em punctured
tree} will be denoted by $t^{\circ}(p)$.
\subsection{Trees for Whitney towers}\label{w-tower-trees-def}
For any Whitney tower $\cW$ define $t(\cW)$ to be the disjoint
union of all the trees $t(p)$ for all unpaired intersection points
$p\in \cW$.
\begin{remn}
As mentioned at the start of this section, $t(\cW)$ can be
enhanced with vertex-orientations and edge-decorations and is
called the {\em geometric intersection tree} of $\cW$. The
invariants $\tau_n$ associated to the obstruction theory mentioned
in Remark~\ref{rem:tau} are determined by such trees (e.g.
\cite{ST3}).
\end{remn}
\subsection{Split Subtowers}
In order to simplify constructions and combinatorics it will be helpful to ``split'' a Whitney
tower into {\em split subtowers} analogous to Krushkal's grope splitting procedure \cite{Kr} which
splits a grope into genus one ``branches'' (higher stages).

\begin{defi}\label{subtower-defi}
A {\em subtower} is a Whitney tower except that the boundaries of the immersed order zero surfaces
in a subtower are allowed to lie in the interior of the 4--manifold. (The boundaries of the order
zero surfaces in a subtower are still required to be embedded.) The notions of {\em order} for
intersection points and Whitney disks are the same as in Definition~\ref{w-tower-defi}.
\end{defi}

\begin{defi}\label{split-subtower-defi}
A subtower $\cW_p$ is {\em split} if it satisfies all of the
following:
\begin{enumerate}
\item
$\cW_p$ contains a single unpaired intersection point $p$,
\item
the order zero surfaces of $\cW_p$ are all embedded 2--disks,
\item
the Whitney disks of $\cW_p$ are all embedded,
\item
the interior of any surface in $\cW_p$ either contains $p$ or
contains a single Whitney arc of a Whitney disk in $\cW_p$,
\item
$\cW_p$ is connected (as a 2--complex in the 4--manifold).

\end{enumerate}
The tree $t(\cW_p)$ associated to such a split subtower $\cW_p$ is
just the tree $t(p)$ determined by its unpaired intersection point
$p$ (\ref{int-point-trees}).
\end{defi}
Note that the order of a split subtower $\cW_p$ is equal to the
order of $p$ and that a normal thickening of $\cW_p$ in the
ambient 4--manifold is just a 4--ball neighborhood of the embedded
tree $t(\cW_p)=t(p)$.
\begin{figure}[ht!]
        \centerline{\includegraphics[scale=.45]{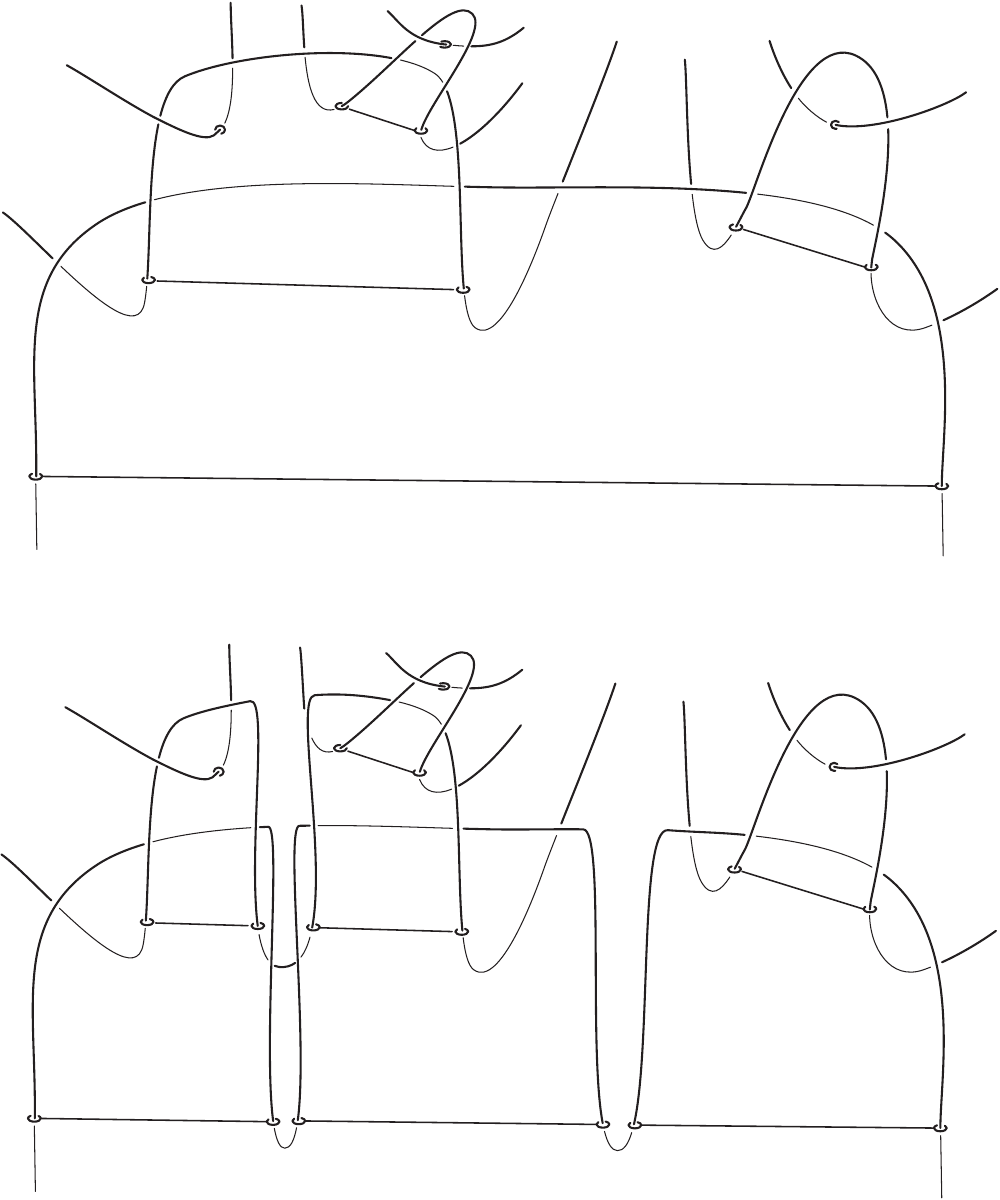}}
        \caption{Part of a Whitney tower before (top) and after (bottom) applying
        the splitting procedure described in the proof of Lemma~\ref{split-tower-lem}.}
        \label{splitting-and-split-tower-fig}

\end{figure}

\subsection{Split Whitney towers and their trees}\label{split-w-tower-tree-subsec}

\begin{lem}\label{split-tower-lem}
Let $\cW$ be a Whitney tower on order 0 surfaces $A_i$. Then, for
any regular neighborhood $N(\cW)$ of $\cW$, there exists a Whitney
tower $\cW_{\mathrm{split}}$ contained in $N(\cW)$ such that:
\begin{enumerate}
\item The order~0 surfaces $A'_i$ of $\cW_{\mathrm{split}}$ only differ from the
$A_i$ by finger moves.
\item All unpaired intersection points of $\cW_{\mathrm{split}}$ are contained
in disjoint split subtowers on sheets of the $A'_i$.
\item $t(\cW)$ and $t(\cW_{\mathrm{split}})$ are isomorphic.
\end {enumerate}
\end{lem}

Such a Whitney tower $\cW_{\mathrm{split}}$ will be called a {\em
split} Whitney tower. The disjoint union of trees
$t(\cW_{\mathrm{split}})$ (\ref{w-tower-trees-def}) associated to
$\cW_{\mathrm{split}}$ is the disjoint union of the trees
$t(\cW_p)$ associated to the split subtowers in
$\cW_{\mathrm{split}}$. Also, $t(\cW_{\mathrm{split}})$ sits as an
{\em embedded} subset of $\cW_{\mathrm{split}}$ via the
description in \ref{int-point-trees}.

\begin{proof}
Starting with the highest order Whitney disks of $\cW$, apply
finger moves as indicated in
Figure~\ref{splitting-and-split-tower-fig}. Working down through
the lower order Whitney disks yields the desired
$\cW_{\mathrm{split}}$.
\end{proof}
This decomposition of a Whitney tower into split subtowers corresponds to the idea that the
disjoint union of trees associated to the Whitney tower captures its essential structure. The next
lemma can be interpreted as justifying that this essential structure is indeed captured by the {\em
un}punctured trees rather than the punctured trees, in the sense that a punctured edge
(corresponding to an unpaired intersection point) can be ``moved'' to any other edge of its tree.
\begin{lem}\label{subtower-lemma}
Let $\cW$ be a split subtower on order zero sheets $s_i$, with
unpaired intersection point $p=W_I\cap W_J\subset \cW$. Denote by
$\nu(\cW)$ a normal thickening of $\cW$, so that $\partial
s_i\subset\partial\nu(\cW)\subset\nu(\cW)\cong B^4$. If $I'$ and
$J'$ are any brackets such that $t(I')\cdot t(J')=t(p)=t(I)\cdot
t(J)$, then after a homotopy (rel $\partial$) of the $s_i$ in
$\nu(\cW)$ the $s_i$ admit a split subtower $\cW'\subset\nu(\cW)$
with single unpaired intersection point $p'=W_{I'}\cap
W_{J'}\subset T'$.
\end{lem}
\begin{figure}[ht!]
        \centerline{\includegraphics[scale=.5]{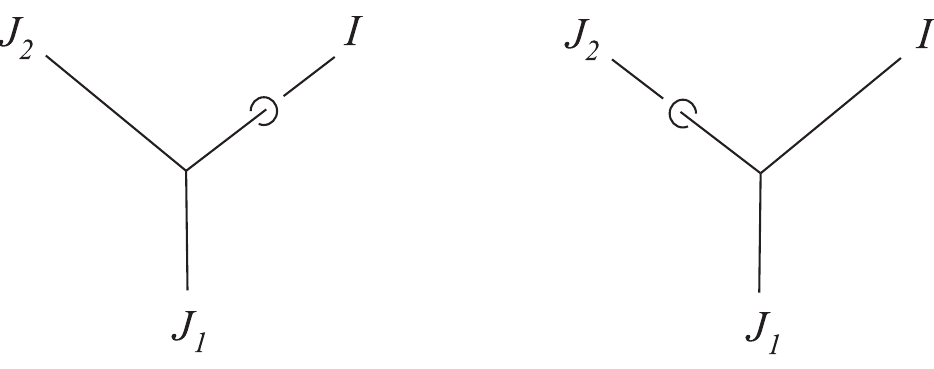}}
        \caption{A local picture of the tree associated to the split
        subtower $\cW$
        before (left) and $\cW'$ after (right) the Whitney move illustrated in
        Figure~\ref{tree-move-and-W-disk-trees-fig}.}
        \label{tree-move-trees-fig}

\end{figure}
\begin{figure}[ht!]
        \centerline{\includegraphics[scale=.5]{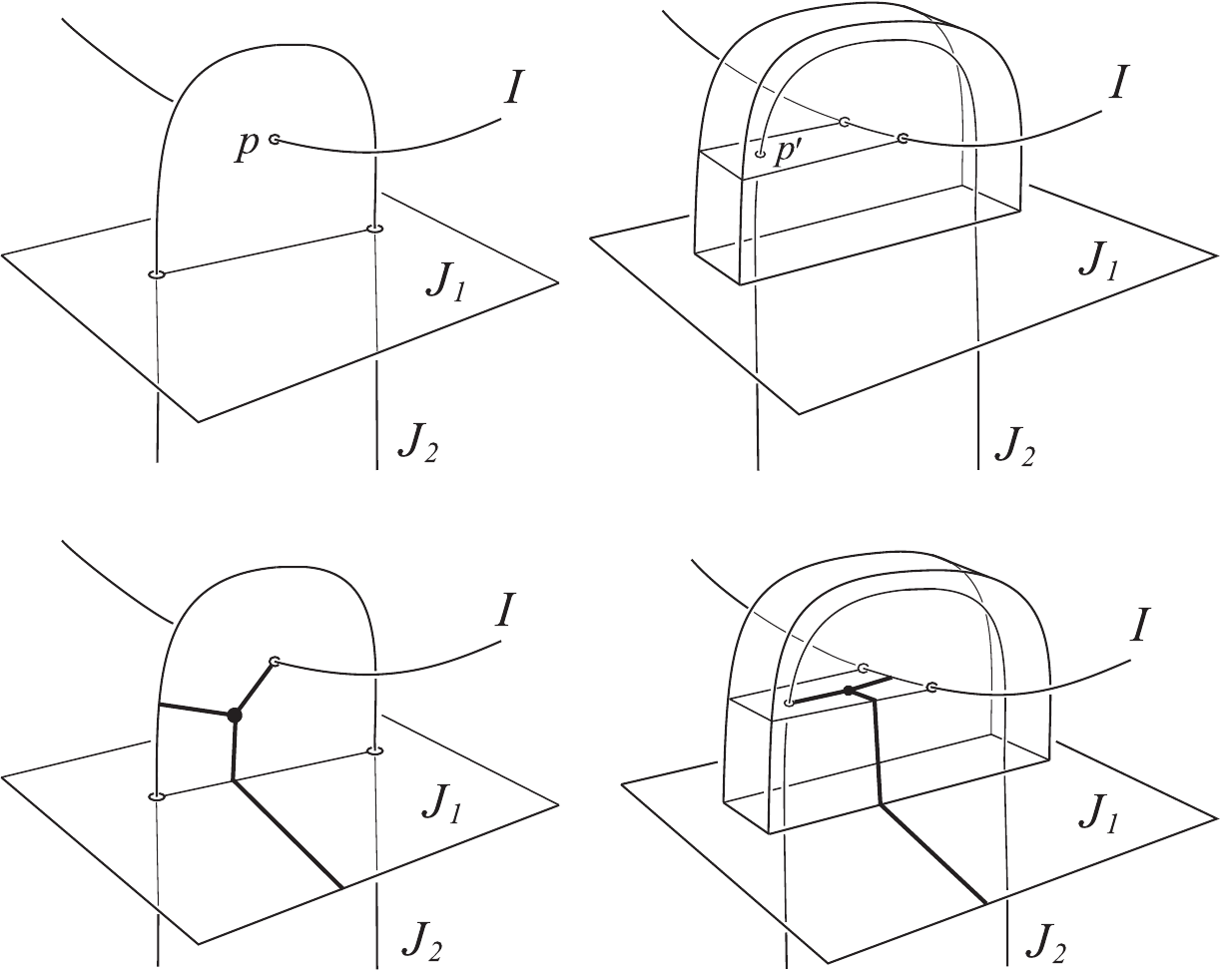}}
        \caption{Upper left, the unpaired intersection point $p=W_I\cap W_J$ in the split
        subtower $\cW$ of Lemma~\ref{subtower-lemma}. Upper right, the unpaired
        intersection point $p'=W_{I'}\cap W_{J'}$ in $\cW'$ after the Whitney move. The lower part of
        the figure shows that the punctured trees differ as indicated in
        Figure~\ref{tree-move-trees-fig}.}
        \label{tree-move-and-W-disk-trees-fig}

\end{figure}
\begin{proof}
It is enough to show that the puncture in $t^{\circ}(p)$ can be ``moved'' to either {\em adjacent}
edge, since by iterating it can be moved to any edge of $t(p)$. Specifically, it is enough to
consider the case where $J=(J_1,J_2)$, $I'=(I,J_1)$ and $J'=J_2$ so that
$I\cdot(J_1,J_2)=(I,J_1)\cdot J_2$ as in Figure~\ref{tree-move-trees-fig}. (Here we are assuming
that $W_J$ is not order zero, since if both $W_I$ and $W_J$ are order zero there is nothing to
prove.) The proof is given by the maneuver illustrated in the upper part of
Figure~\ref{tree-move-and-W-disk-trees-fig}: Use the Whitney disk $W_J$ to guide a Whitney move on
$W_{J_1}$. This eliminates the intersections between $W_{J_1}$ and $W_{J_2}$ (as well as
eliminating $W_J$ and $p$) at the cost of creating a new cancelling pair of intersections between
$W_{J_1}$ and $W_I$. This new cancelling pair can be paired by a Whitney disk $W_{(I,J_1)}$ having
a single intersection point $p'$ with $W_{J_2}$. That this achieves the desired effect on the
punctured tree is clear from the lower part of Figure~\ref{tree-move-and-W-disk-trees-fig}.
\end{proof}

\section{Grope subtowers}\label{sec:grope-subtowers}
In this section a hybrid grope--Whitney tower combination is introduced which will be used to
interpolate between gropes and Whitney towers in the proofs of the next section: A {\em split grope
subtower} is a collection of capped gropes whose caps support certain split subtowers. We assume
from now on that all gropes are dyadic (\ref{dyadic-gropes-defi}).

\begin{figure}[h!]
        \centerline{\includegraphics[scale=.5]{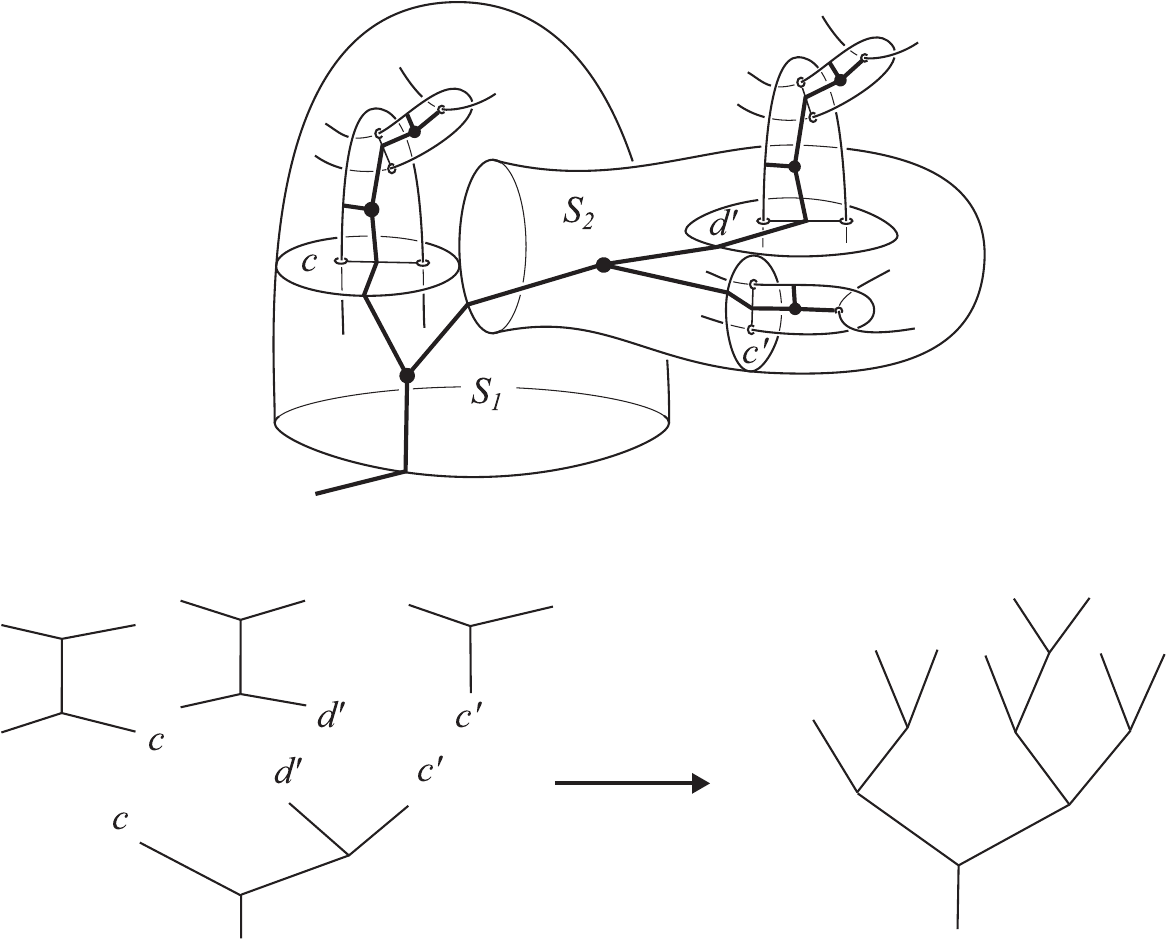}}
        \caption{The tree $t(g^\cW)$ associated to a split grope subtower
        $g^\cW$ is formed from the trees associated to the caps and gropes
        by gluing univalent vertices associated to common caps (bottom)
        and embeds in $g^\cW$ (top).}
        \label{grope-subtower-tree1-with-trees-fig}

\end{figure}

The content of this section can essentially be grasped by inspecting
Figure~\ref{grope-subtower1-fig} and Figure~\ref{grope-subtower-tree1-with-trees-fig}, which
illustrate split grope subtowers and their trees, and observing that the limiting cases reduce to
gropes (\ref{Order 0 subtowers on caps}) and Whitney towers (\ref{Class 1 grope subtowers})
respectively.

\subsection{Split grope subtowers.}\label{grope-subtowers}
\begin{defi}

Let $g_i^c$ be a collection of (dyadic) $A_i$--like capped gropes properly immersed in a
4--manifold such that:
\begin{enumerate}

\item The higher (greater than zero) stages of the $g_i$ are all
disjointly embedded and disjoint from the interiors of all caps,

\item the interiors of all caps are disjointly embedded,

\item each cap $c$ supports a split subtower $\cW_c$ whose other order 0 surfaces are sheets of the
0th stage surfaces $A^0_i$,

\item the $\cW_c$ are disjoint and contain all singularities among the $0$th stages $A^0_i$.

\end{enumerate}
Denote by $g_i^\cW$ the union of the grope $g_i^c$ and the
subtowers on its caps. The union $g^\cW$ of all the $g_i^\cW$ is a
{\em split grope subtower} (see Figure~\ref{grope-subtower1-fig}).
\begin{figure}[ht!]
        \centerline{\includegraphics[scale=.5]{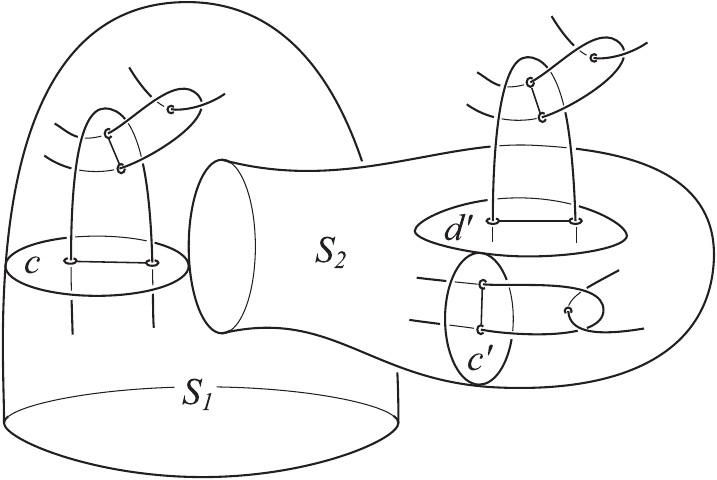}}
        \caption{A {\em split grope subtower} is a collection of properly immersed
        capped gropes,
        with each cap supporting
        a split subtower whose other order 0 surfaces are sheets of the
        0th stages of the gropes.}
        \label{grope-subtower1-fig}

\end{figure}
The {\em class} of $g_i^\cW$ is the class of the underlying grope
$g_i^c$. The {\em class} of $g^\cW$ is the minimum of the classes
of the $g_i^\cW$.

The {\em order} of $g_i^\cW$ is defined inductively as follows: If
$g_i^\cW$ is class 1, then the order of $g_i^\cW$ is the minimum
of the orders of the split subtowers on the caps of $g_i^\cW$ (the
immersed disks that fill in the punctures of $A^0_i$). If
$g_i^\cW$ has class 2 and a single genus one first stage, then the
order of $g_i^\cW$ is the sum of the orders of the split subtowers
on the dual pair of caps of $g_i^\cW$. If $g_i^\cW$ has class 2
and more than one first stage, then the order of $g_i^\cW$ is the
minimum of the orders of the first stages. If the class of
$g_i^\cW$ is greater than 2, then the order of $g_i^\cW$ is
defined (inductively) to be the minimum of the sums of the orders
of the pairs of dual grope subtowers that are attached to the
first stages of $g_i^\cW$. The {\em order} of $g^\cW$ is the
minimum of the orders of the $g_i^\cW$.

\end{defi}

\subsection{Trees for dyadic split grope subtowers.}
For each $g_i^\cW$ in a split grope subtower $g^\cW$, construct
the disjoint union of (rooted) unitrivalent trees $t(g_i^\cW)$
from $t(g_i^c)$ (defined in \ref{grope-trees}) by gluing on the
trees $t(\cW_c)$ (defined in \ref{subtower-defi}) along the
univalent vertices that correspond to caps. Specifically, a
univalent vertex of $t(g_i^c)$ which corresponds to a cap $c$ in
$g_i^c$ is identified with the univalent vertex of $t(\cW_c)$
which corresponds to $c$, where $\cW_c$ is the subtower on $c$.
This identification is to a single non-vertex point in an edge of
$t(g_i^\cW)$ (see the lower part of
Figure~\ref{grope-subtower-tree1-with-trees-fig}). Doing this for
all caps on $g_i$ and all $i$ yields all the $t(g_i^\cW)$. The
disjoint union of trees $t(g^\cW)$ associated to the split grope
subtower $g^\cW$ is defined to be the disjoint union of the
$t(g_i^\cW)$ and sits as a subset of $g^\cW$ (upper part of
Figure~\ref{grope-subtower-tree1-with-trees-fig}).

\subsection{Order 0 split grope subtowers.}\label{Order 0 subtowers on
caps} If a class $m$ split grope subtower $g^\cW$ has order 0,
then all the split subtowers $\cW_c$ in $g^\cW$ are order 0 which
just means that each cap of every $g_i^c$ has exactly one interior
intersection point with a sheet of some $A^0_j$. In this case, the
trees $t(g_i^\cW)$ and $t(g_i^c)$ are clearly isomorphic for all
$i$ and each univalent vertex corresponds to a sheet of some
$A^0_j$.

\subsection{Class 1 split grope subtowers.}\label{Class 1 grope
subtowers}

If each $g_i^\cW$ in an order $m$ grope subtower $g^\cW$ has class
1, then the caps fill in the punctures in the $0$th stages $A^0_i$
to form the order 0 surfaces in an order $m$ split Whitney tower
$\cW$ on immersions of the $A_i$ extending the embedded $A^0_i$.
The disjoint unions of trees $t(g^\cW)$ and $t(\cW)$ are
isomorphic, with the root of each chord in $t(g_i^c)$
corresponding to an $i$-labelled vertex of a tree in $t(\cW)$.

\section{Proof of Theorem~\ref{class-order-thm}}\label{sec:thms}
The equivalence of the statements in Theorem~\ref{class-order-thm}
in the introduction follows directly from the more detailed
Theorems \ref{tower-to-grope-thm} and \ref{grope-to-tower-thm}
which are stated and proved in this section. A key element of
these theorems is that when passing between gropes and Whitney
towers, the associated trees are ``preserved''. In this setting,
an {\em isomorphism} between rooted and unrooted (disjoint unions
of) trees will always mean an isomorphism between the underlying
unrooted trees, but will also include a correspondence between the
roots and certain specified univalent vertices, e.g. the roots in
$t(g_i)$ will always correspond to $i$-labelled univalent vertices
of $t(\cW)$ when passing between gropes $g_i$ and a Whitney tower
$\cW$ on order zero surfaces $A_i$. (These isomorphisms also
preserve the {\em signed} trees associated to gropes and Whitney
towers as in \cite{CST,WTCCL,ST3}.)
\subsection{From Whitney towers to
gropes.}\label{w-tower-to-grope-subsec}
\begin{thm}\label{tower-to-grope-thm}
Let $\cW$ be an order $(n-1)$ Whitney tower on properly immersed
surfaces $A_i$ in a 4--manifold $X$. Then, for any regular
neighborhood $N(\cW)$ of $\cW$, there exist class $n$ $A_i$--like
capped gropes $g_i^c$ in $X$ such that:
\begin{enumerate}
\item The $0$th stage $A_i^0$ of each $g_i^c$ is $A_i$ minus (perhaps) some
sheets containing Whitney arcs or intersection points of $\cW$ in
$A_i$.

\item The union of the $g_i^c$ are contained in $N(\cW)$.

\item The $g_i^c$ have disjoint properly embedded bodies $g_i$.

\item Each cap of every $g_i^c$ has a single interior intersection
with some $A^0_j$.

\item $t(\cW)$ is isomorphic to the disjoint union of the
$t(g_i^c)$, with $j$-labelled univalent vertices in $t(\cW)$
corresponding to either vertices in the $t(g_i^c)$ associated to
caps which intersect $A_j^0$ or roots in $t(g_j^c)$; furthermore,
it may be arranged that this isomorphism takes any chosen
preferred $i$-labelled univalent vertices on the trees in $t(\cW)$
to the root vertices of the trees in $t(g_i^c)$.
\end{enumerate}
\end{thm}
The proof of Theorem~\ref{tower-to-grope-thm} is well illustrated by
Figure~\ref{tubing-W-disk1and2with-trees-fig} together with the observation that the pictured case
can always be arranged by Lemma~\ref{subtower-lemma}.
\begin{proof}
First split $\cW$ (Lemma~\ref{split-tower-lem}) so that $t(\cW)$
is the disjoint union of the split subtower trees $t(\cW_p)$ each
of order at least $n-1$. For each $\cW_p$, choose a preferred
univalent vertex of $t(\cW_p)$ and let $A^0_i$ denote the
punctured surfaces which are the complements of the sheets of the
$A_i$ that correspond to the chosen preferred vertices. (Each of
these sheets is either a neighborhood of a Whitney disk boundary
arc or a neighborhood of an unpaired intersection point.) These
chosen vertices will end up corresponding to root vertices in the
$t(g_i^c)$ (which are associated to the 0th stages of the capped
gropes $g_i^c$) as in statement (v) of the theorem.

Now $\cW$ is a grope subtower $g^\cW$ of class 1 and order
$(n-1)$: The $0$th stages of $g^\cW$ are the $A^0_i$ and the caps
of $g^\cW$ are the sheets of the $A_i$ that correspond to the
chosen preferred vertices. The trees $t(\cW)$ and $t(g^\cW)$ are
isomorphic. In particular, the theorem is true for $n=1$ by
\ref{Order 0 subtowers on caps}.

Assume now that $n\geq 2$. The proof will be completed by the
following construction which shows how to decrease the order of
$g^\cW$ while increasing the class of $g^\cW$ in a manner that
preserves the tree $t(g^\cW)$. When each cap supports an order 0
split subtower the proof is done by \ref{Order 0 subtowers on
caps}.
\begin{figure}[ht!]
        \centerline{\includegraphics[scale=.5]{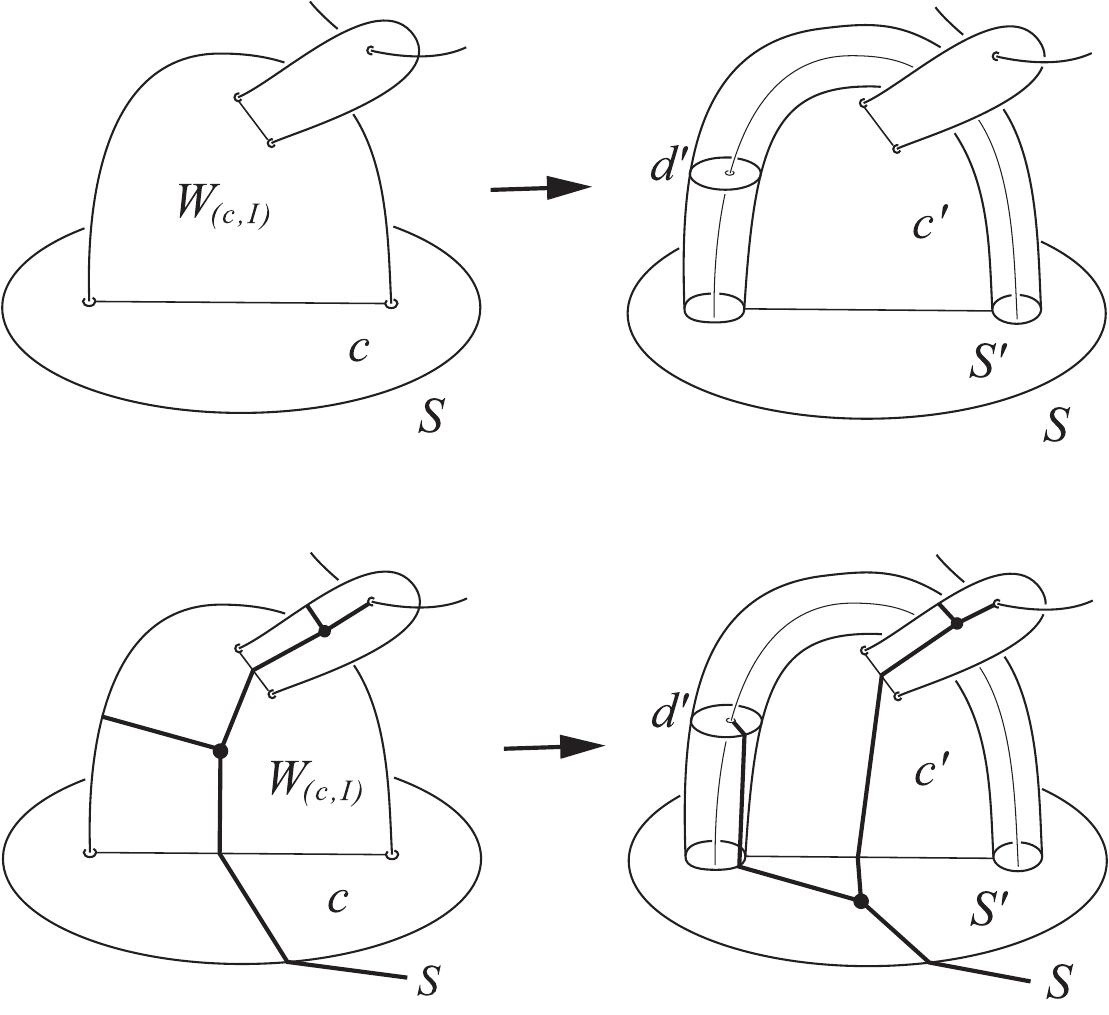}}
        \caption{Upper left to right: Reducing the order and increasing the class of
        a split grope subtower $g^\cW$ by 0-surgering a cap $c$ along the boundary arc of a Whitney disk.
        That $c$ contains such a boundary arc (rather than the unpaired intersection
        point of $\cW_c$) can be arranged by Lemma~\ref{subtower-lemma}.
        Lower left to right: This procedure preserves $t(g^\cW)$.}
        \label{tubing-W-disk1and2with-trees-fig}

\end{figure}
\subsection{Decreasing the order and increasing the class and of a
split grope subtower.} \label{decreasing-order-subsec} Consider a
cap $c$ attached to some stage $S$ in a grope subtower $g^\cW$
such that the order of the split subtower $\cW_c$ supported by $c$
is greater than or equal to 1. There are two cases to consider:
Either $c$ contains a boundary arc of a Whitney disk in $\cW_c$,
or $c$ contains the unpaired intersection point $p$ of $\cW_c$.

First assume that $c$ contains a boundary arc of a Whitney disk
$W_{(c,I)}$ (see upper left of
Figure~\ref{tubing-W-disk1and2with-trees-fig}). In this case,
``tube'' (0-surger) $c$ along the (other) boundary arc of
$W_{(c,I)}$ that lies in $W_I$ as indicated in the upper right of
Figure~\ref{tubing-W-disk1and2with-trees-fig}. This changes $c$
into a genus one capped surface stage $S'$. One cap $c'$ is
$W_{(c,I)}$ minus a small collar and the other dual cap $d'$ is a
meridional disk to $W_I$. Both of these caps support split
subtowers of order strictly less than $\cW_c$ since the trees
$t(\cW_{c'})$ and $t(\cW_{d'})$ are gotten from $t(\cW_c)$ by
removing the edge adjacent to the vertex associated to $c$ and
cutting $t(\cW_c)$ at the vertex associated to $W_{(c,I)}$. The
tree associated to the new grope subtower is the same as the
original tree $t(g^\cW)$ since the effect of creating $S'$ from
$c$ just isotopes (in $X$) the trivalent vertex (basepoint) of
$t(g^\cW)$ in $W_{(c,I)}$ down to a trivalent vertex (basepoint)
in $S'$ as illustrated in the lower part of
Figure~\ref{tubing-W-disk1and2with-trees-fig}.

Now assume that $c$ contains the unpaired intersection point $p$
of $\cW_c$. We may assume that $p$ is the intersection between $c$
and a Whitney disk $W_{(I,J)}$ since $\cW_c$ has order greater
than or equal to 1. Modify $\cW_c$ by one iteration of the
procedure of Lemma~\ref{subtower-lemma}: Do the $W_{(I,J)}$
Whitney move on $W_I$. This creates a cancelling pair of
intersections between $c$ and $W_I$ which are paired by a Whitney
disk $W_{(c,I)}$ that has a single intersection with $W_J$ (as in
Figure~\ref{tree-move-and-W-disk-trees-fig} but with different
labels). The modified split subtower $\cW'_c$ on $c$ has the same
tree as $\cW_c$ and we are back to the previous case where $c$
contains a boundary arc of $W_{(c,I)}$. (We remark that this step
where $c$ contains the unpaired intersection point $p$ of $\cW_c$
could alternatively be handled by building higher grope stages out
of Clifford tori as in the {\em grope duality} constructions of
\cite{KrT}).)
\end{proof}

\subsection{From gropes to Whitney
towers.}\label{grope-to-w-tower-subsec} The gropes in the statement of
Theorem~\ref{class-order-thm} can be arranged to satisfy the hypotheses of next theorem by using
finger moves to push down cap-intersections (2.5 of \cite{FQ}) into the order 0 surfaces and by
applying Krushkal's grope splitting technique \cite{Kr}.

\begin{thm}\label{grope-to-tower-thm}
Let $g_i^c$ be a collection of class $n$ $A_i$--like dyadic capped
gropes in a 4--manifold $X$ with disjoint properly embedded bodies
$g_i$ such that all the caps have disjointly embedded interiors
and each cap contains only a single interior intersection point
with some $0$th stage $A_j^0$. Then, for any regular neighborhood
$N(g_i^c)$ of the union of the $g_i^c$, there exists an order
$(n-1)$ Whitney tower $\cW$ in $X$ such that:
\begin{enumerate}

\item $\cW$ is contained in $N(g_i^c)$.

\item The order 0 surfaces of $\cW$ are immersions of the $A_i$
extending the embeddings $A_i^0$ up to regular homotopy (rel $\partial$).

\item $t(\cW)$ is isomorphic to the disjoint union of the
$t(g_i^c)$ with $j$-labelled univalent vertices in $t(\cW)$
corresponding to either vertices in the $t(g_i^c)$ associated to
caps which intersect $A_j^0$ or roots in $t(g_j^c)$.

\end{enumerate}
\end{thm}
The proof of Theorem~\ref{grope-to-tower-thm} is well illustrated
by Figure~\ref{surgered-cap-1and2-with-trees-fig} together with
the observation that the pictured case can always be arranged by
Lemma~\ref{subtower-lemma}:
\begin{figure}[ht!]
        \centerline{\includegraphics[scale=.5]{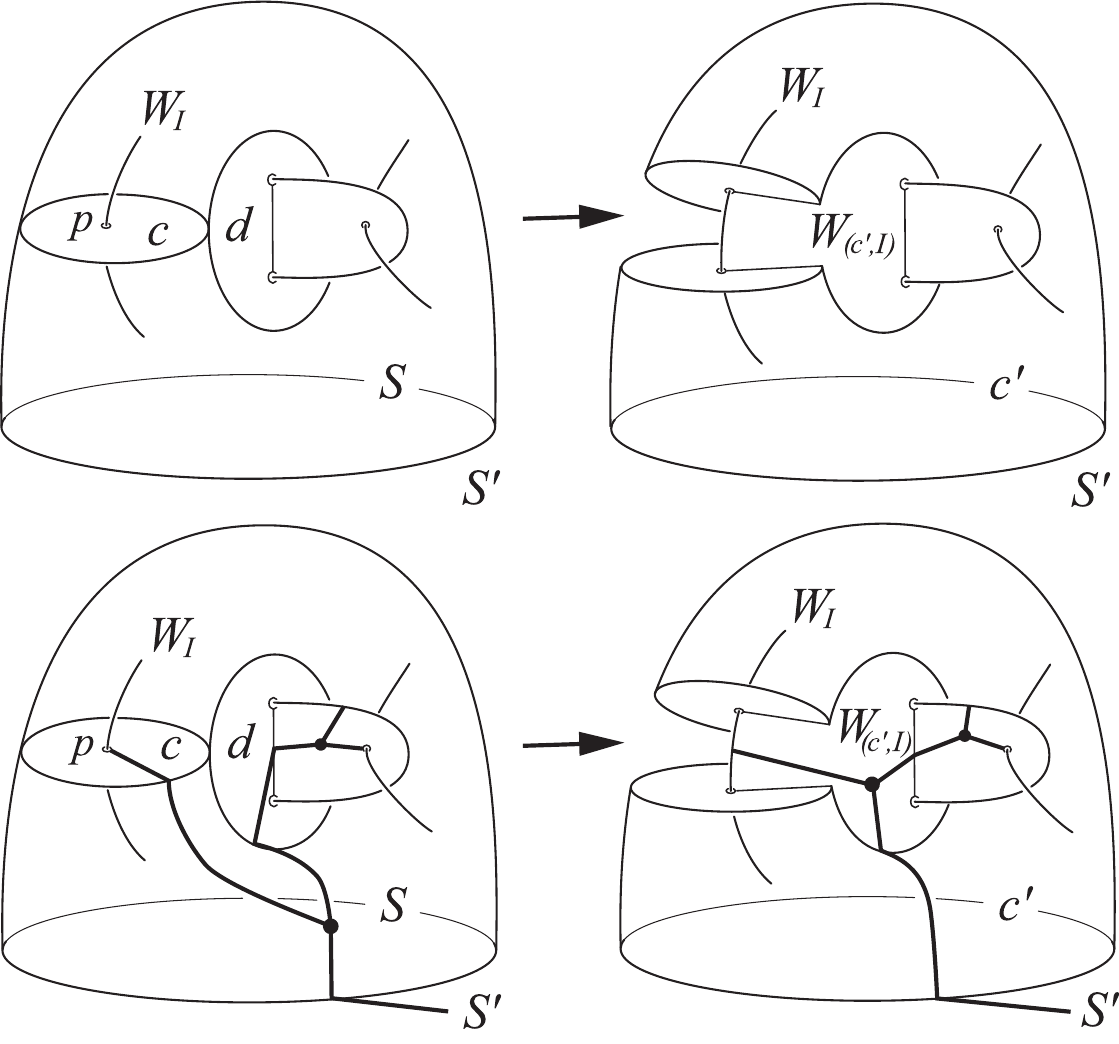}}
        \caption{Upper left to right: Surgering one of a dual pair of
        caps in a split grope subtower $g^\cW$ reduces the class and increases the order of
        $g^\cW$. This procedure requires the surgered cap $c$ to contain the unpaired intersection point
        $p$ in its split subtower $\cW_c$ as can be arranged by
        Lemma~\ref{subtower-lemma}. Lower left to right: This
        surgery preserves $t(g^\cW)$.}
        \label{surgered-cap-1and2-with-trees-fig}

\end{figure}
\begin{proof}
When $n=1$, the $g_i^c$ form a grope subtower $g^\cW$ of class 1
and the Theorem is true by \ref{Class 1 grope subtowers}.

Assuming $n\geq 2$, the proof is completed by the following
construction (essentially the inverse to the construction
\ref{decreasing-order-subsec} in the proof of
Theorem~\ref{tower-to-grope-thm}) which decreases the class and
increases the order of a grope subtower $g^\cW$ while preserving
the associated trees. When each $g_i^c$ in $g^\cW$ has class 1 the
proof is complete by \ref{Class 1 grope subtowers}.

\subsection{Decreasing the class and increasing the order of a
grope subtower.} Let $c$ and $d$ be a pair of dual caps on a
surface stage $S$ in a grope subtower $g^\cW$ supporting split
subtowers $\cW_c$ and $\cW_d$. By applying
Lemma~\ref{subtower-lemma} to $\cW_c$, we may arrange that $c$
contains the unpaired intersection point $p=c\cap W_I$ in $\cW_c$.
Using $c$ to ambiently surger $S$ changes $S$ into a cap $c'$ on
the stage $S'$ below $S$. This new cap $c'$ has a cancelling pair
of intersections with $W_I$ (due to $p$, the intersection $c$ had
with $W_I$), which can be paired with a Whitney disk $W_{(c',I)}$
formed from the old cap $d$ by attaching a thin band as pictured
in the upper right of
Figure~\ref{surgered-cap-1and2-with-trees-fig}. The cap $c'$
supports a split subtower $\cW_{c'}$ and, as illustrated in the
lower part of Figure~\ref{surgered-cap-1and2-with-trees-fig},
there is no change in $t(g^\cW)$ since the effect of the surgery
is just to isotope (in $X$) the unitrivalent vertex of $t(g^\cW)$
that was in $S$ up into $W_{(c',I)}$. (Here we are still denoting
the modified grope subtower by $g^\cW$.) Repeated application of
this construction eventually eliminates all dual pairs of caps so
that each $g_i^c$ in $g^\cW$ has class 1.
\end{proof}


\section{Proof of
Corollary~\ref{cor:height}}\label{sec:height-cor} The idea of the
proof of Corollary~\ref{cor:height} in the introduction is to use
Theorem~\ref{grope-to-tower-thm} of
Subsection~\ref{grope-to-w-tower-subsec} to convert the symmetric
grope into a Whitney tower whose associated trees are all
symmetric and then use Lemma~\ref{subtower-lemma} to create the
desired Whitney tower by moving the unpaired intersection points
appropriately. The gropes in this section are disk-like and
assumed to be dyadic as can always be arranged by Krushkal's
splitting procedure \cite{Kr}.

\subsection{Symmetric gropes and trees}\label{sym-grope-tree-subsec}
\begin{figure}[ht!]
        \centerline{\includegraphics[scale=.5]{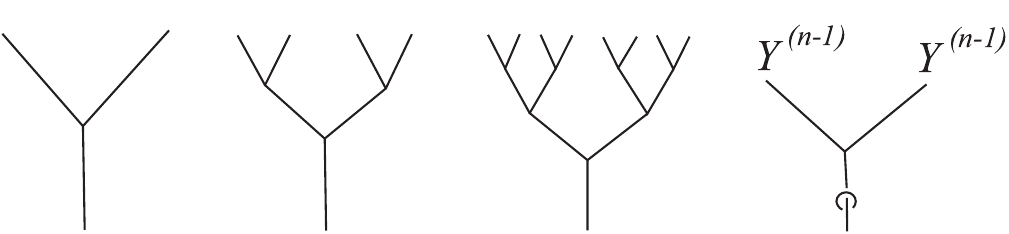}}
        \caption{From left to right: The trees $Y^1$, $Y^2$, $Y^3$ and
        (a punctured) $Y^n$.}
        \label{Y-trees-fig}

\end{figure}
Let $Y^1$ denote the order $1$ rooted $Y$--shaped tree that
corresponds to a punctured torus viewed as a grope of class 2
(with the bottom univalent vertex designated as the root). For
integers $n > 1$, define $Y^{(n+1)}$ to be the rooted product
$Y^n\ast Y^n$ (see Figure~\ref{Y-trees-fig}). A grope $g$ is {\em
symmetric} and has {\em height} $n$ if all the trees in $t(g)$ are
of the form $Y^n$. For $n\geq 1$, define the tree $Y^{(n.5)}$ to
be the rooted product $Y^{(n-1)}\ast Y^n$ (where $Y^0$ is the
rooted chord). A grope $g$ has {\em height} $n.5$ if all the trees
in $t(g)$ are of the form $Y^{(n.5)}$. Note that a grope of height
$n$ (resp. $n.5$) has class $2^n$ (resp.
$2^n+2^{(n-1)}=(1.5)(2^n)$).

\subsection{The height of a Whitney
tower}\label{w-tower-height-def} Translating the definition of
height given in \cite{COT} into our language we have: An
order~$(2^n-2)$ Whitney tower $\cW$ has {\em height} $n$ ($n\geq
1$) if the interiors of all Whitney disks in $\cW$ only intersect
surfaces of the same order. Thus, the lowest order unpaired
intersections in a Whitney tower of height $n$ are of order
$(2^{(n-1)}-1)+(2^{(n-1)}-1)=2^n-2$ and occur among the highest
order Whitney disks (of order $2^{(n-1)}-1$).

A Whitney tower of {\em height} $n.5$ ($n\geq 1$) is a Whitney
tower of height $n$, together with order $2^n-1$ Whitney disks
pairing all order $2^n-2$ intersection points with the requirement
that the interiors of these order $2^n-1$ Whitney disks may
intersect each other and the order $2^{(n-1)}-1$ Whitney disks but
are disjoint from all surfaces of lower order. A Whitney tower of
height $n.5$ has order $2^n+2^{(n-1)}-2$.

\begin{proof}(of Corollary~\ref{cor:height}) Applying
Theorem~\ref{grope-to-tower-thm} of Section~\ref{sec:thms} to a
disk-like grope $g$ of height $n$ (resp. $n.5$) yields a Whitney
tower $\cW$ of order $2^n-1$ (resp. $2^{(n-1)}+2^n-1$) with
$t(\cW)$ isomorphic to $t(g)$ so that all the connected trees in
$t(\cW)$ are of the form $Y^n$ (resp. $Y^{(n.5)}$). After
splitting $\cW$ (Lemma~\ref{split-tower-lem}) and applying
Lemma~\ref{subtower-lemma}, we may arrange that the punctured edge
in each punctured tree $t^{\circ}(p)$ in $t^{\circ} (\cW)$ is
adjacent to (what was) a root vertex, that is, the only unpaired
intersection points of $\cW$ occur between the order zero 2--disk
and Whitney disks whose associated trees are of the form $Y^n$
(resp. $Y^{(n.5)}$) as in the far right of
Figure~\ref{Y-trees-fig}. The Whitney disks of $\cW$ correspond to
the trivalent vertices of $t(\cW)$ and one can check by examining
the shape of the $Y^n$ trees which make up $t^{\circ} (\cW)$ that
$\cW$ satisfies the above definition of height $n$ (resp. $n.5$).
In fact, $\cW$ satisfies the stronger condition that its
intersections of order $2^n-2$ (between Whitney disks of order
$2^{(n-1)}-1$), which are allowed to be unpaired in \cite{COT},
are in fact all paired by order $2^n-1$ Whitney disks (each of
which corresponds to the trivalent vertex adjacent to the root of
a $Y^n$). This is as expected by Theorem~\ref{class-order-thm}
since $\cW$ should have order $2^n-1$. However, these order
$2^n-1$ Whitney disks intersect the order zero 2--disk so that
$\cW$ does {\em not} have height $n+1$. The case of half-integer
height $n.5$ is checked similarly.
\end{proof}


\section{Proof of
Corollary~\ref{cor:half-grope} and the Whitney move IHX construction}\label{sec:half-grope-cor}
This section contains a proof of Corollary~\ref{cor:half-grope} which is based on a geometric
realization of the IHX Jacobi relation in the setting of Whitney towers (Lemma~\ref{IHX-lemma}
below).
\begin{figure}[ht!]
        \centerline{\includegraphics[scale=.6]{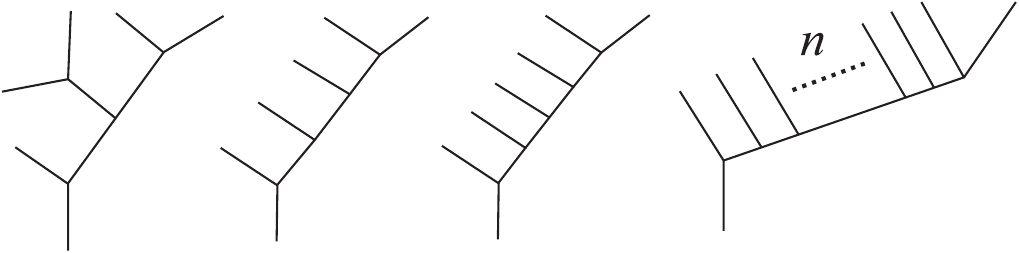}}
        \caption{From left to right: The non-simple tree of lowest order
        (order 4) and the simple trees of order 4, 5,
        and $6+n$.}
        \label{simple-trees-fig}

\end{figure}
A dyadic (capped) $A$--like grope $g$ is a {\em half-grope} if all the trees in $t(g)$ are {\em
simple} (right- or left-normed) as illustrated in Figure~\ref{simple-trees-fig}; note that the
roots (shown pointing down in the figure) are required to be at an ``end'' of the tree.
\begin{proof} 
To prove Corollary~\ref{cor:half-grope}, choose caps for the class
$n$ $g_i$ (which we may assume are dyadic) and use
Theorem~\ref{grope-to-tower-thm} of Section~\ref{sec:thms} to
convert the $g_i^c$ into a Whitney tower $\cW$ of order $n-1$. By
the following Proposition~\ref{simple-tower-prop}, $\cW$ can be
modified (rel boundary) to a Whitney tower $\cW'$ (of the same
order) with $t(\cW')$ consisting of only simple trees. Then
converting $\cW'$ back into a grope via
Theorem~\ref{tower-to-grope-thm} of Section~\ref{sec:thms} yields
the desired half gropes $h_i$ in $X$ bounded by $\gamma_i$. (Here
we are using (v) of Theorem~\ref{tower-to-grope-thm} to send
chosen univalent ``end vertices'' to roots, while preserving
trees.)
\end{proof}

\begin{prop}\label{simple-tower-prop}
Let $\cW$ be any order $n$ Whitney tower on order 0 surfaces
$A_i$. Then, after a regular homotopy (rel $\partial$), the $A_i$
admit an order $n$ {\em simple} Whitney tower $\cW'$ contained in
a neighborhood of $\cW$, that is, $t(\cW')$ consists of only
simple trees.
\end{prop}
The proof of Proposition~\ref{simple-tower-prop} uses the geometric IHX Lemma~\ref{IHX-lemma} below
to follow the algebraic proof that the usual group of unitrivalent trees occurring in finite type
theory is spanned by simple trees as given in e.g.~\cite {B2,CT1}:
\begin{proof}The simple trees are characterized by the property that any
maximal length chain of edges contains every trivalent vertex. Let
$t(p)=t(\cW_p)\in t(\cW)$ be a tree associated to a split subtower
$\cW_p\subset \cW$ contained in a Whitney tower $\cW$ (which we may
assume is split by Lemma~\ref{split-tower-lem}). If $t(p)$
contains a trivalent vertex $v_1$ which is of distance 1 away from
a trivalent vertex $v_0$ contained in some maximal chain of edges
not containing $v_1$, then we may assume, by
Lemma~\ref{subtower-lemma}, that $v_0$ corresponds to $W_{(I,J)}$
and $v_1$ corresponds to $W_{((I,J),K)}$ which is incident to some
other sheet $W_L$ in $\cW_p$ as in left hand side of
Figure~\ref{IHX-trees-fig}. Lemma~\ref{IHX-lemma} below shows how
to modify $\cW$ near $\cW_p$ so as to replace $t(p)\in t(\cW)$ by
the two trees on the righthand side of Figure~\ref{IHX-trees-fig}
having the same order as $t(p)$ but with longer length edge
chains. By iterating this modification (for all components of
$t(\cW)$) we eventually arrive at the desired $\cW'$ with all
components of $t(\cW')$ simple trees.
\end{proof}
\begin{figure}[ht!]
        \centerline{\includegraphics[scale=.5]{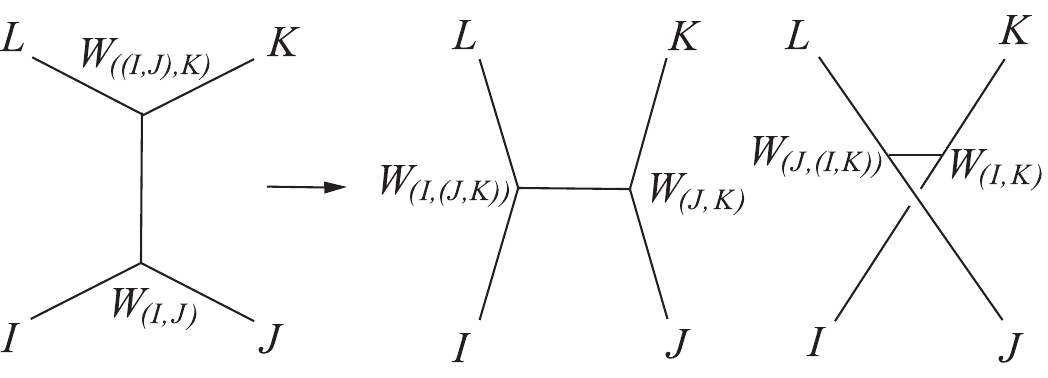}}
        \caption{The IHX relation for Whitney disks in a split subtower
        replaces a split subtower whose tree looks locally like the one on the left
        with a pair of nearby disjoint split subtowers whose trees
        look locally like the trees on the right.}
        \label{IHX-trees-fig}

\end{figure}
\begin{lem}[Geometric IHX--Whitney move version]\label{IHX-lemma}
Let $\cW_p$ be a split subtower in a split Whitney tower $\cW$.
Let $W_{((I,J),K)}$ be a Whitney disk in $\cW_p$ so that
$t(\cW_p)$ looks locally like the leftmost tree in
Figure~\ref{IHX-trees-fig}. Then $\cW$ can be modified in a
regular neighborhood $\nu(\cW_p)$ of $\cW_p$ yielding a split
Whitney tower $\cW'$, on the same order 0 surfaces, with $\cW_p$
replaced by disjoint split subtowers $\cW_{p'}$ and $\cW_{p''}$
contained in $\nu(\cW_p)$ such that the trees $t(\cW_{p'})$ and
$t(\cW_{p''})$ are as pictured on the right hand side of
Figure~\ref{IHX-trees-fig}.
\end{lem}
\begin{figure}[ht!]
        \centerline{\includegraphics[scale=.65]{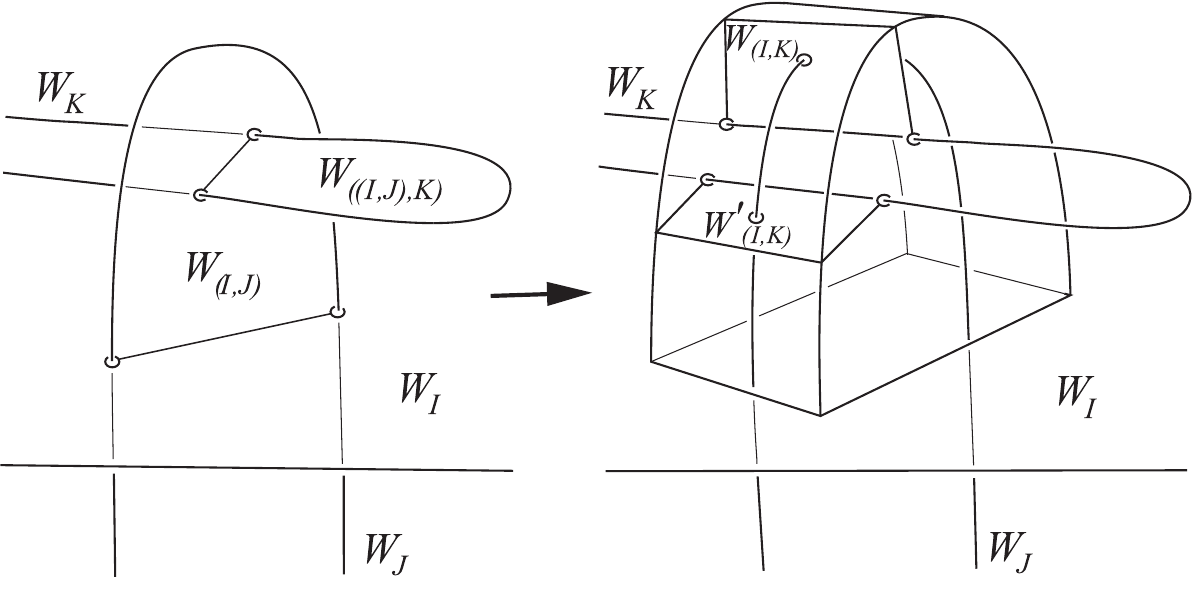}}
        \caption{The IHX construction starts with a $W_{(I,J)}$ Whitney move
        on $W_I$. Note that intersections between $W_{((I,J),K)}$ and $W_L$ are {\em not} shown in this
        figure (and are suppressed in subsequent figures as well).}
        \label{IHX-W-move1and2-fig}

\end{figure}
The modification involves Whitney moves, finger moves and taking
parallel copies of some of the Whitney disks in $\cW_p$. The
reader familiar with the orientation and sign conventions of
\cite{ST3} can check by inserting signs and orientations in the
figures that the following construction actually replaces an ``I''
tree with the difference ``$\mathrm{H}-\mathrm{X}$'' as in the
usual IHX relation of finite type theory. Note that this differs
from the closely related 4--dimensional IHX construction in
\cite{CST} which creates the trees
$\mathrm{I}-\mathrm{H}+\mathrm{X}$ for a Whitney tower on
2--spheres in 4--space by modifying the boundaries of Whitney
disks.
\begin{figure}[ht!]
        \centerline{\includegraphics[scale=.55]{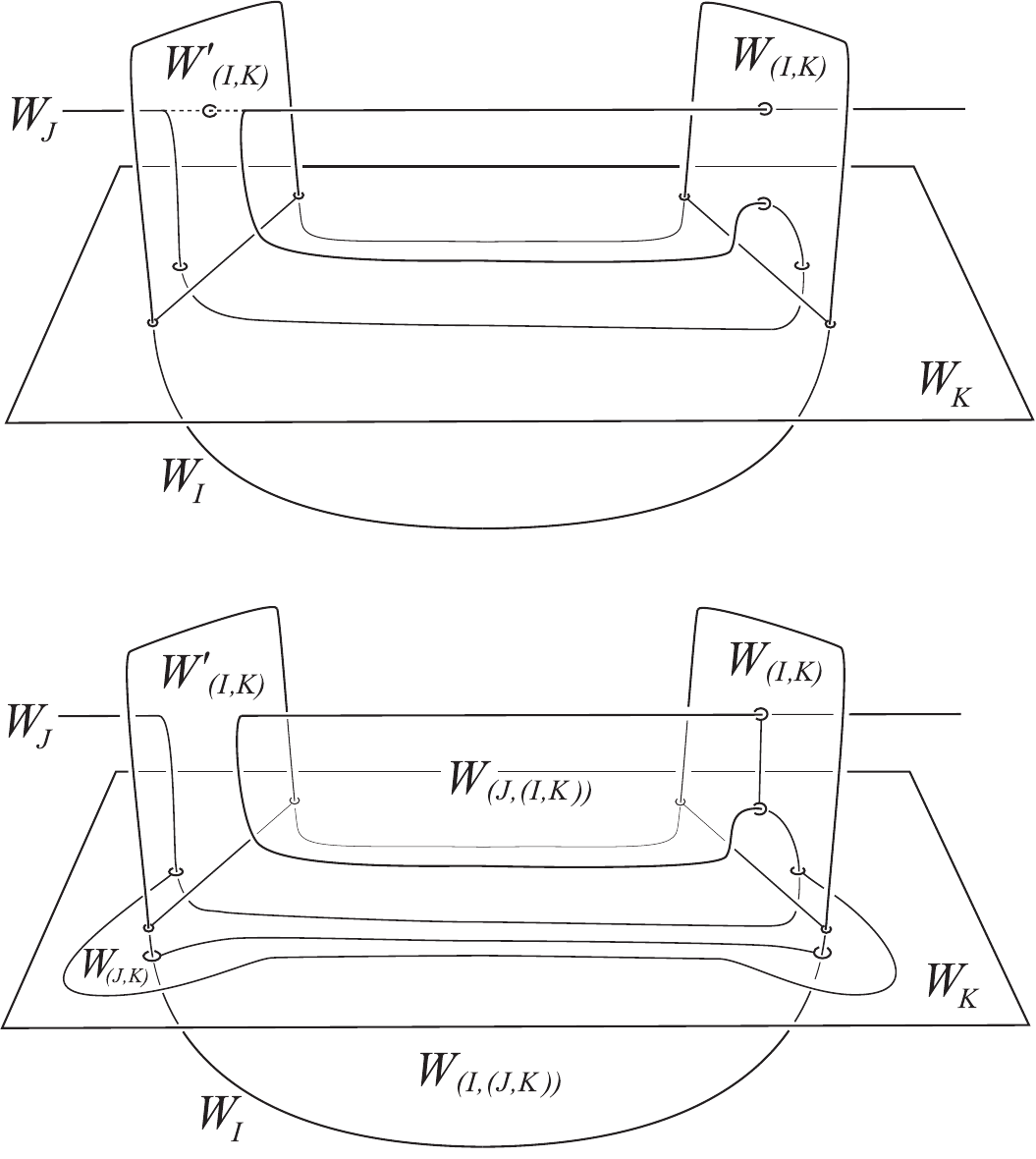}}
        \caption{The intersection point $W_J\cap W_{(I,K)}'$ is
        `transferred' via a finger move (top)
        to create a cancelling pair $W_J\cap
W_{(I,K)}$ paired by $W_{(J,(I,K))}$ at the cost of also creating $W_J\cap W_K$ paired by
$W_{(J,K)}$ and $W_I\cap W_{(J,K)}$ paired by $W_{(I,(J,K))}$ (bottom).}
        \label{IHX-transfer-move1and1A-fig}

\end{figure}
\begin{figure}[ht!]
        \centerline{\includegraphics[scale=.55]{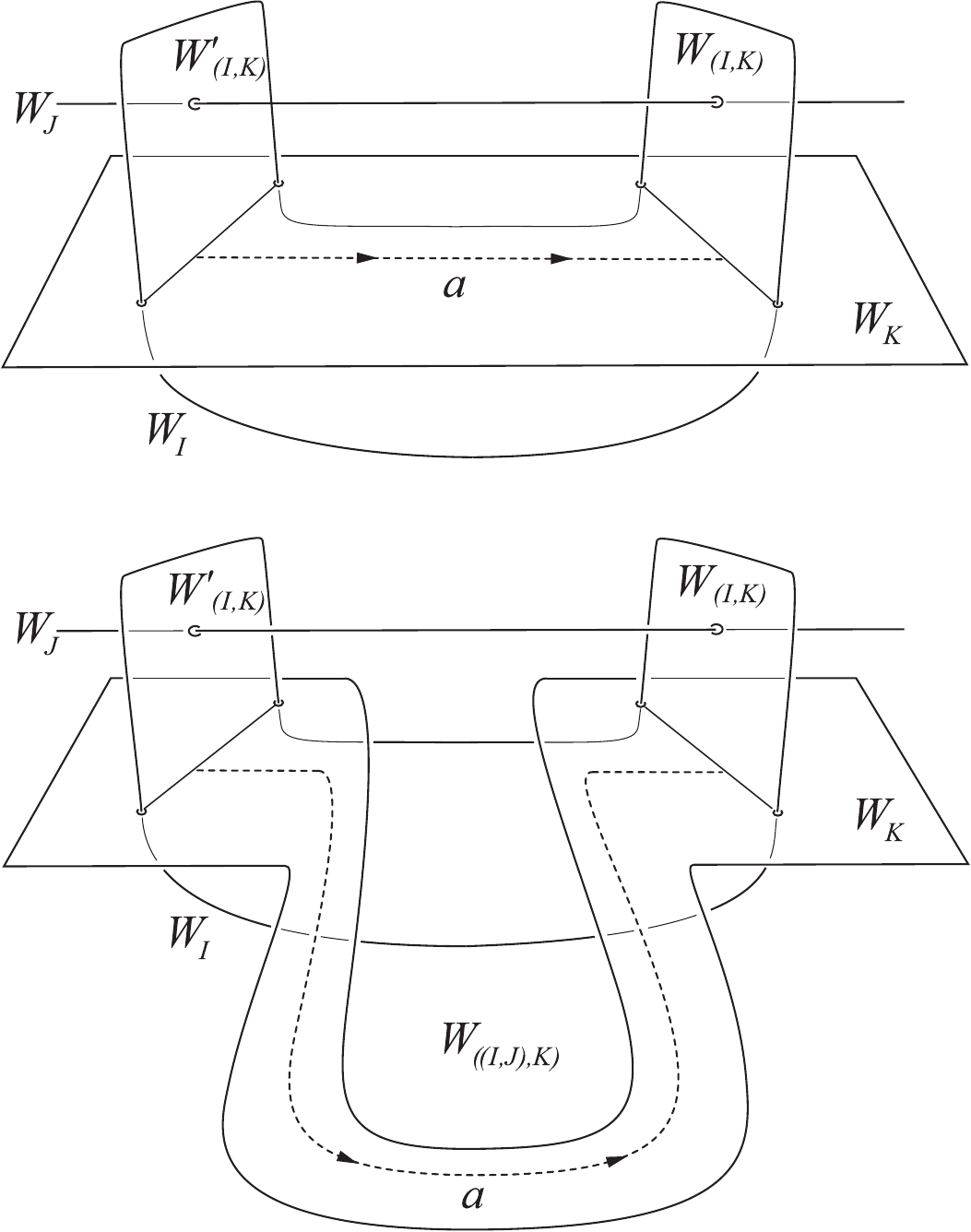}}
        \caption{The transferring finger move is guided by an arc $a$ (top)
        which can be taken to run along what used to be the part of the boundary arc of
        $W_{((I,J),K)}$ lying in $W_K$ (bottom). Indicated in the bottom picture is where
         $W_{((I,J),K)}$ used to be.}
        \label{IHX-transfer-move2and3-fig}

\end{figure}
\begin{proof}
The first step in the modification is to do the $W_{(I,J)}$
Whitney move on $W_I$ (see Figure~\ref{IHX-W-move1and2-fig}) and
disregard, for the moment, the Whitney disks in the part of $\cW$
corresponding to the sub-tree $L$. This eliminates the cancelling
pair of intersections between $W_I$ and $W_J$ at the cost of
creating two cancelling pairs of intersections between $W_I$ and
$W_K$ which we pair by Whitney disks $W_{(I,K)}$ and $W_{(I,K)}'$
as illustrated in Figure~\ref{IHX-W-move1and2-fig}. The new
Whitney disks $W_{(I,K)}$ and $W_{(I,K)}'$ each have a single
interior intersection with $W_J$ and the next step is to
``transfer'' (as illustrated in the upper part of
Figure~\ref{IHX-transfer-move1and1A-fig}) the intersection point
$W_J\cap W_{(I,K)}'$ to create a cancelling pair $W_J\cap
W_{(I,K)}$ paired by $W_{(J,(I,K))}$ at the cost of also creating
$W_J\cap W_K$ paired by $W_{(J,K)}$ and $W_I\cap W_{(J,K)}$ paired
by $W_{(I,(J,K))}$ (as illustrated in the lower part of
Figure~\ref{IHX-transfer-move1and1A-fig}). Note that
Figure~\ref{IHX-transfer-move1and1A-fig} differs from
Figure~\ref{IHX-W-move1and2-fig} by a rotation of coordinates
which brings the sheet of $W_K$ into the ``present'' slice of
3--space.
\begin{figure}[ht!]
        \centerline{\includegraphics[scale=.55]{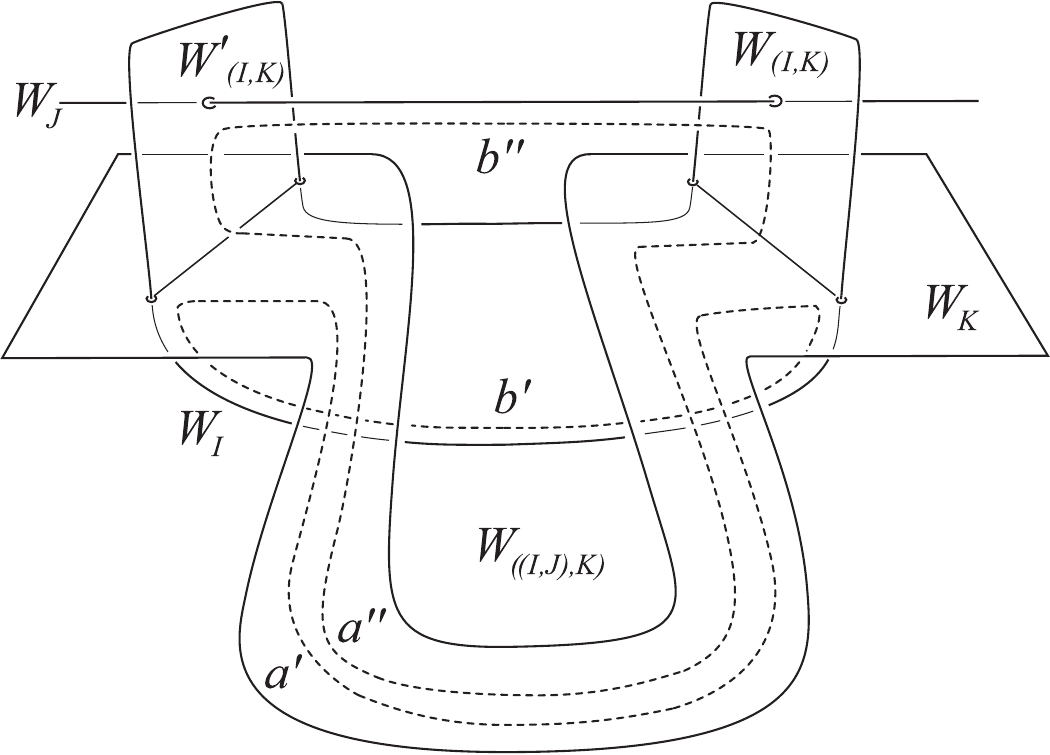}}
        \caption{Before the transfer move: New Whitney disks $W_{(I,(J,K))}$
        and $W_{(J,(I,K))}$, whose boundaries are the unions of arcs $a'\cup b'$
        and $a''\cup b''$ (see also Figure
        ~\ref{IHX-W-move1and2B-fig}),
         will be created from parallel copies of the old $W_{((I,J),K)}$.}
        \label{IHX-transfer-move3B-fig}
\end{figure}
This transfer move was described in \cite{Y} (see also \cite{ST1,ST3}) 
and is just a (non-generic) finger move applied to
$W_J$. The important thing to note here is that the finger move is
guided by an arc $a$ (see Figure~\ref{IHX-transfer-move2and3-fig})
from $\partial W_{(I,K)}'$ to $\partial W_{(I,K)}$ in $W_K$ and we
can take this arc to run along what used to be the part of
$\partial W_{((I,J),K)}$ lying in $W_K$. This is illustrated in
the lower part of Figure~\ref{IHX-transfer-move2and3-fig} which
gives a better picture of the situation before the finger move is
applied.
\begin{figure}[ht!]
        \centerline{\includegraphics[scale=.55]{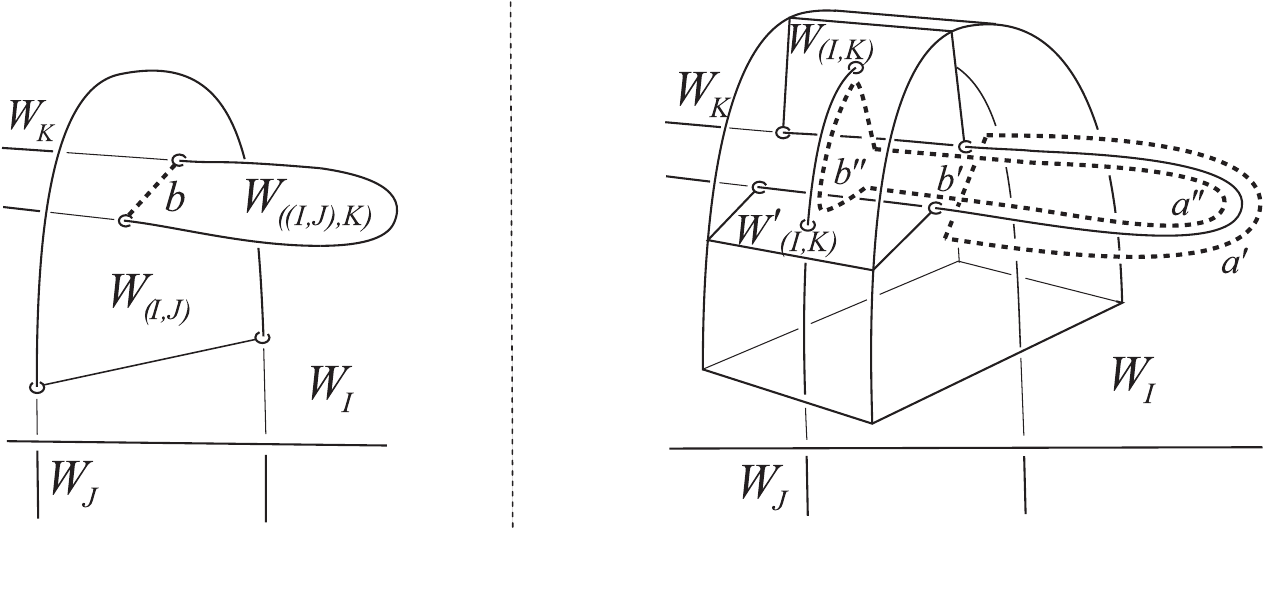}}
        \caption{Applying the transfer move to the right-hand side will create
        new Whitney disks $W_{(I,(J,K))}$
        and $W_{(J,(I,K))}$, whose boundaries are the unions of arcs $a'\cup b'$
        and $a''\cup b''$ (see also Figure
        ~\ref{IHX-transfer-move3B-fig}),
         from parallel copies of the old $W_{((I,J),K)}$ shown on the left.}
        \label{IHX-W-move1and2B-fig}

\end{figure}
The Whitney disks $W_{(I,(J,K))}$ and $W_{(J,(I,K))}$ can be taken
to be parallel copies of the old $W_{((I,J),K)}$ as follows: The
boundary of $W_{(I,(J,K))}$ (resp. $W_{(J,(I,K))}$) consists of
arcs $a'$ and $b'$ (resp. $a''$ and $b''$) where $a'$ and $a''$
are tangential push-offs of $a$ in $W_K$ and $b'$ and $b''$ are
normal push-offs of what was the boundary arc $b$ of
$W_{((I,J),K)}$ in $W_{(I,J)}$. This is shown in both
Figure~\ref{IHX-transfer-move3B-fig} and
Figure~\ref{IHX-W-move1and2B-fig}, where again it is easier to
picture things before the transferring finger move. Since
$W_{((I,J),K)}$ was framed and embedded, $W_{(I,(J,K))}$ and
$W_{(J,(I,K))}$ can be formed from two disjoint parallel copies of
$W_{((I,J),K)}$ which each intersect $W_L$ as $W_{((I,J),K)}$ did.
Two parallel copies of each Whitney disk that was in the part of
$\cW_p$ corresponding to $L$ can be used to recover the order of
the original Whitney tower (by pairing the new intersections
``over'' $W_{(I,(J,K))}$ and $W_{(J,(I,K))}$ corresponding to
$L$), which means that exactly two new unpaired intersection
points $p'$ and $p''$ have been created with corresponding trees
$t(p')$ and $t(p'')$ as shown locally in the right hand side of
Figure~\ref{IHX-trees-fig}. After the transferring finger move,
the $W'_{(I,K)}$ Whitney move can be done (on either sheet)
without affecting anything else. Finally, $W_I$, $W_J$ and $W_K$
will need to be split since they now each contain two boundary
arcs of Whitney disks. Splitting $W_I$, $W_J$ and $W_K$ down into
the lower order Whitney disks (as in Lemma~\ref{split-tower-lem})
yields the two split subtowers $\cW_{p'}$ and $\cW_{p''}$.

\end{proof}


\section{Proof of Corollary~\ref{cor:k-slice}\label{sec:k-slice}}
In this section the main theorems together with the geometric IHX construction of the previous
section are used to prove Corollary~\ref{cor:k-slice} in the introduction. We refer the reader to
\cite{C1} or \cite{IO} for the formal definition of {\em $k$-null-cobordism} or {\em $k$-slice}. It
is enough to show that the link components in $S^3=\partial B^4$ bound disjointly embedded surfaces
$A_i$ such that each $A_i$ contains a symplectic basis of circles which bound continuous maps of
class $k$ gropes in $B^4\setminus \bigcup_iA_i$. In fact, we will find such $A_i$ with symplectic
bases of circles bounding {\em embedded} class $k$ gropes in $B^4\setminus \bigcup_iA_i$. The same
proof shows that class $2k$ grope concordant links are $k$-cobordant.
\begin{figure}[h]
        \centerline{\includegraphics[scale=.6]{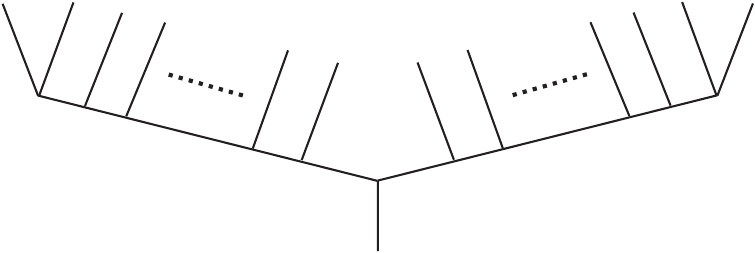}}
        \caption{A simple tree associated to an order $2k-1$ (or greater)
        intersection point in an order $2k-1$ Whitney tower $\cW$ with
        a preferred univalent vertex
        (at least) $k-1$ trivalent vertices away from both ends. Using
        Theorem~\ref{tower-to-grope-thm} of Section~\ref{sec:thms}
        to convert $\cW$ to disjoint gropes
        $g_i$ with all such preferred univalent vertices going to roots in $t(g_i)$
        yields $k$-slicing surfaces (the bottom stages of the $g_i$).}
        \label{k-slice-tree-fig}
\end{figure}
\begin{proof}
Apply Krushkal's grope-splitting procedure \cite{Kr} and
Theorem~\ref{grope-to-tower-thm} of Section~\ref{sec:thms} to the
class $2k$ (disk-like) gropes (as in the hypotheses of
Corollary~\ref{cor:k-slice}) to get an order $2k-1$ Whitney tower
on immersed 2--disks bounded by the link components in $B^4$. Then
apply Proposition~\ref{simple-tower-prop} to get a Whitney tower
$\cW$ with $t(\cW)$ consisting of only simple trees. Since the
order of $\cW$ is $2k-1$, the order of each tree in $t(\cW)$ is at
least $2k-1$. This means that we can specify a preferred univalent
vertex in each simple tree in $t(\cW)$ such that the trivalent
vertex adjacent to the preferred univalent vertex is at least
$k-1$ trivalent vertices away from both ends of the simple tree,
as illustrated in Figure~\ref{k-slice-tree-fig}. Now use
Theorem~\ref{tower-to-grope-thm} of Section~\ref{sec:thms} to
convert back to disjointly embedded (disk-like) gropes $g_i$ with
the preferred univalent vertices in $t(\cW)$ going to the root
vertices in the $t(g_i)$. As is seen in the shape of the trees
$t(g_i)$ (Figure~\ref{k-slice-tree-fig}), the desired $A_i$ are
the class 2 sub-gropes formed by the bottom (0th and 1st) stages
of the $g_i$ and the higher ($\geq 2$) stages form class $k$
gropes which are attached to symplectic bases on the $A_i$ and
disjointly embedded in $B^4\setminus \bigcup_iA_i$.
\end{proof}


\end{document}